\pdfoutput=1

\documentclass[12pt,a4paper]{article}

\usepackage[english]{babel}
\usepackage{geometry}    
\usepackage{amsmath,amssymb,amsfonts}
\usepackage[round]{natbib}
\usepackage{graphicx,xcolor}

\usepackage{overpic}
\usepackage{bbm}

\usepackage{vmargin}
\setpapersize{A4}
\setmarginsrb{2cm}{2cm}{2cm}{2cm}{2pt}{8pt}{0pt}{1cm}

\DeclareMathOperator*{\argmin}{arg min}

\newcommand{\inv}{^{-1} }

\usepackage{algorithm,algorithmic}

\usepackage{authblk}

\title{Truncation map estimation based on bivariate probabilities 
and validation for the truncated plurigaussian model}

\author[1,2]{Astrakova, Alina} 
\author[1]{Oliver, Dean S.}
\author[3]{Lantu\'{e}joul, Christian}
\affil[1]{Uni Research (Center for Integrated Petroleum Research), Bergen, Norway\\ alina.astrakova@uni.no}
\affil[2]{University of Bergen, Bergen, Norway} 
\affil[3]{Mines ParisTech, Fontainebleau, France}

\begin{document}

\maketitle


\vspace{1cm}

\begin{abstract}

The truncated plurigaussian model is often used to simulate the spatial distribution of random categorical variables such as geological facies. The problems addressed in this paper are the estimation of parameters of the truncation map for the truncated plurigaussian model. Unlike standard truncation maps, in this paper a colored Voronoi tessellation with number of nodes, locations of nodes, and category associated with each node all treated as unknowns in the optimization. Parameters were adjusted to match categorical bivariate unit-lag probabilities, which were obtained from a larger pattern joint distribution estimates from the Bayesian maximum-entropy approach conditioned to the unit-lag probabilities. The distribution of categorical variables generated from the  estimated truncation map  was close to the target unit-lag bivariate probabilities.

The predictive performance of the model is evaluated using scoring rules, and conditioning of the latent Gaussian fields to log-data is generalized for the case when the truncated bigaussian model is governed by a colored Voronoi tessellation of the truncation map. 

\end{abstract}

\section{Introduction}

Initially developed by \citet{matheron:87}, the categorical-valued truncated Gaussian model provides a method for simulating random categorical-valued fields with desired proportions and approximate transition probabilities 
from  a latent Gaussian random field (GRF). The generality  of the model was increased when it was 
extended to the truncated plurigaussian (TPG) model \citep{galli:94,leloch:94} in which simulation of a categorical variable is based on the values of an arbitrary number of GRFs. In fact, however, most plurigaussian applications should more accurately be called bigaussian, as in most applications a categorical variable is assigned based on values of a random standard Gaussian pair. The function used to assign the categorical value is commonly rather simple so that it can be  visualized by coloring the practical domain of $\mathbb{R}^2$ with one color per category. The hyperparameters for the function that maps the GRFs to a categorical variable and the function itself are  referred to as the truncation map. To simulate a categorical-valued random field in a bigaussian case, two GRFs are sampled. The spatial correlation within each of the GRF and possible correlation between those fields represent other model parameters.

The categorical-valued TPG model  has seen increasing application in a number of fields, mostly due to Gaussian distribution of the latent random fields which is useful in several data-assimilation applications. 
The attractiveness of the TPG model is also partly due to the ability to handle large model-size and partly due to the advances in modeling real geological fields and real data. 
One area with significant  TPG application is in the characterization of heterogeneous aquifers. \citet{cherubini:09} were able to reproduce complex categorical lithofacies geometry using the truncated bigaussian model for a contaminated aquifer in which category proportions varied with depth. A similar TPG aquifer characterization problem was described by \citet{mariethoz:09}  with porosity heterogeneity simulated as a GRF for each category. \citet{perulero:12,perulero:14} assessed the  suitability of the TPG model through improved aquifer flow simulation responses.  
Apart from aquifers, \citet{emery:10} provided  a model for mineral proportion evaluation in an ore deposit.  
Another extensive area of TPG application is related to petroleum reservoir characterization. 
Fault facies modeling with TPG was provided by \citet{fachri:13}. \citet{armstrong:11} include several complex reservoir examples of primary diagenesis effects characteristic of carbonate sedimentary systems. \citet{carrillat:10} compared the use of TPG on a giant carbonate oilfield with application of other methods. \citet{albertao:05} and \citet{alanezi:13} used TPG to describe the distribution of facies for complex reservoir models while using varying category proportions computed from seismic data.

The application of the TPG method to unconditional simulation (or conditional to a few observations) of categorical fields is often fairly straightforward once the model parameters have been estimated, but estimation of the model parameters can be difficult. Although it is clear that the 
TPG parameters should be estimated with respect to available static data (logs and cores), there is no general way to do this, and those data  are typically not directly involved in the parameter optimization. If the observations are abundant, their derivatives, for example, experimental variograms, cross-variograms and proportions, can be computed. Otherwise, if the observations are sparse or are not representative, some categorical space relations might be estimated from other sources, often coming from geological studies of analogue fields. 
A simplified approach to estimation of the truncation map is to chose a truncation map  from some benchmark set \citep[among others]{armstrong:11,galli:06}. The choice of map from the truncation set is based mostly on allowing or not-allowing some category contacts, while classical thresholds of the truncation map have either vertical or horizontal orientation with a given mutual arrangement.
The exact threshold values are then computed from category proportions. An advantage of this method of the truncation map parametrization is that it is trivial to  adjust locally varying category proportions while the other parameters are fixed.
To handle more complex contacts between categorical variables, \citet{xu:06} increased the  number of GRFs but retained the classical thresholding style.
\citet{allard:12} developed a more general approach to truncation map estimation using  kernel regression methods to relate categorical variables to the auxiliary latent variables. The truncation map in their method is non-parametric, allowing complex contacts between categorical variables, but the method requires joint observation of the categorical and the Gaussian variables at the well locations. In another approach which avoided the assumption of rectangular regions in the truncation map, \citet{deutsch:14} attempted optimization of the truncation map parametrized with Voronoi tessellation. One node was used for each category while optimization of node locations was split into a series of steps that operated on subsets of the data, for example, proportions and fixed-lag bivariate  probabilities.

In addition to estimating parameters of the truncation map, it may be necessary to estimate the parameters of the GRF such as the covariances of the latent fields and the possible correlation between them. These parameters must be estimated jointly with parameters of the truncation map, hence the problem of estimating 
the Gaussian covariance matrices is almost necessarily iterative. 
Covariance estimation usually assumes a fixed covariance model (known in advance), stationarity and geometric anisotropy. In the approach by \citet{xu:06} it was possible to make an iterative adjustment for variogram ranges and angles. \citet{kyriakidis:99} studied an accurate way to estimate Gaussian variograms in case of a univariate model with two categories separated by one threshold, based on experimental categorical variograms as input.

A key requirement of the TPG method is to be able to condition the GRF to actual categorical observations. This is necessary  either because of the need to simulate categorical fields conditional to well observations, or to estimate the parameters of the truncation map conditional to the data.
When the parameters of the TPG model are established, the conditioning of the Gaussian variables related to the categorical observations can be considered as a separate problem that is commonly solved with Gibbs sampler. For a large number of observations, however, the sampler encounters a problem related to moving search neighborhood application and covariance matrix inversion. 
This problem  was addressed for the unconditional random fields by \citet{lantuejoul:12} using a propagative version of the Gibbs sampler. 
The approach was adapted to the TPG fields with linear inequality constraints by \citet{emery:14}.  

The problem that is addressed in this paper is the need to reproduce geocellular models by 
multiple samples from the joint distribution of the categorical values (facies) to which petrophysical properties can be conditioned. 
The ensemble of sampled models can then be used to quantify the uncertainty in reservoir flow behavior.
The method is based on the following assumptions. The prior model parameter distribution for the truncation map is
based on a parametrization using colored Voronoi tessellations. 
This parametrization is more flexible than the classical truncation map parametrization, which is mostly based on proportions and on allowing (or not) facies contacts.
The prior data parameter distribution ( following the notation by \citet{tarantola:05}) is based on 
the Bayesian maximum-entropy (BME) estimate of a small joint distribution of few observations (a pattern), used for regularization of the data, which itself consists of unit-lag bivariate probabilities.
The goal is  to estimate the parameters of a TPG model (or models) that provides categorical samples/patterns that are indistinguishable from the BME estimates.
The comparison of the two models pattern distribution is made on the basis of expected frequencies from the BME estimate to observed frequencies  in realizations from the truncated multivariate model.

First, the truncation map parametrized as  colored Voronoi tessellation is estimated  while assuming that parameters of the GRF variograms are known. 
Second, to handle the problem of conditioning correlated Gaussian random variables to categorical observations when the relationship is governed by the Voronoi truncation map,  we generalize the constrained version of the propagative Gibbs sampler by \citet{emery:14}. 
If a truncation map has complex form or isolated areas of the same category, the standard Gibbs sampler might not converge to a likely state.
The alternative sampler can more easily condition the correlated observation to the categories presented as unions of  disconnected polygonal areas in the truncation map.
The method proved to converge rapidly and allows faster sampling of multiple Gaussian random vectors that all give the correct conditioning.
Third, the validation of the parameter estimates is based on the scoring rule performance as regards to a materialized event of joint categorical observations, mimicking possible geological data of consequent geological facies observations at well locations.

\section{Methods}
\subsection{Methods: Truncation model estimation}

\subsubsection{Prior information on model parameters}
\label{sec:pr_info_mod_par}

The distribution of the categorical random vector is governed by the truncated plurigaussian model. 
For $C$ being a finite set of facies categories, a truncation map with parameters $\theta$ is a map 
$ M_{\theta}: \mathbb{R}^{2} \mapsto C$.
The set of parameters $\theta$ is given below for a particular type of parametrization. 
$M_\theta$ maps each bivariate Gaussian realization $(x,y) \in \mathbb{R}^{2}$ to the set $C$.
When Gaussian vector realizations,  
\begin{equation}
\label{eq:gaus_vec}
x_A = (x_\alpha, \alpha \in A) ,
y_A = (y_\alpha, \alpha \in A)
\end{equation}
defined on a vector space $A \subset \mathbb{R}^{2}$ (or $\mathbb{R}^{3})$, are used instead, the map $M_{\theta}(x_A,y_A)= z_A$ represents a categorical vector realization $z$, obtained from component pair mapping,
$M_{\theta}(x_\alpha,y_\alpha)= z_\alpha, \alpha \in A$.
Vector $z_A = (z_\alpha, \alpha \in A)$ is a random realization of  
a random vector $Z_A = (Z_\alpha, \, \alpha \in A)$, and $x_A,y_A$ are realizations of the so-called latent variables $X_A,Y_A$.

Although the standard parametrization of a truncation map for a truncated bigaussian model
is via threshold values of the latent Gaussian variables, the truncation map can alternatively
be parametrized by a colored Voronoi tessellation \citep{du:99}. 
The prior information for the map $M_\theta$ is conveniently specified in terms of a distribution for the number of nodes of the tessellation, $\mathcal{T}$, node coordinates $(\chi_\tau,\upsilon_\tau) \in \mathbb{R}^{2}$ and their categories $\zeta_{\tau}$. 
We assume that the number of nodes $\mathcal{T}$ follows the Poisson law $P_\mu(\cdot)$, with specified mean of the Poisson distribution, $\mu$ and that the location of
each coordinate $\chi_\tau,\upsilon_\tau, \tau =1,\dots, \mathcal{T}$, has a prior standard normal distribution, while node categories $\zeta_{\tau}$ are independent and uniformly distributed over the set of categories $C$. 
It is assumed that $C$ includes all possible categories of the field $C$. 
Thus, the set of parameters $\theta$ includes all the above mentioned
\begin{equation}
\theta = (\chi_\tau,\upsilon_\tau, \zeta_\tau, \ \tau =1,\dots, \mathcal{T}).
\label{eq:model_set_par}
\end{equation}
Its prior distribution density take the following form.
\begin{equation}
f_\Theta(\theta) = \prod\limits_{\tau =1}\limits^\mathcal{T}  g(\chi_\tau) g(\upsilon_\tau) \frac{1}{|C|},
\label{eq:prior_mod}
\end{equation}
where $g(\cdot)$ is the standard normal distribution density; and where $|C|$ denotes the cardinality of the set $C$.
The mapping of each pair of latent Gaussian values $(x_\alpha,y_\alpha)$ to a categorical variable using the Voronoi tessellation of the truncation map is done as
\begin{equation}
M_{\theta}(x_\alpha,y_\alpha) = \zeta_{\tau_*}, \; \tau_*= \argmin_\tau \| (x_\alpha,y_\alpha) - (\chi_\tau,\upsilon_\tau) \|.
\end{equation}
In other words, the mapping assigns the category at location $\alpha$ to be equal to the category of the node closest to the latent Gaussian pair related to this location.
The covariance matrices of the latent GRFs $X_A,Y_A$ are assumed to be known. Different hypotheses about the covariance matrices structures could possibly be compared with the validation techniques discussed below. However, this work only shows the results for truncation map validation for specified covariance matrices.

In this work $X_A,Y_A$ are assumed to be independent, but the estimation method equally could have  been applied to the case where the two Gaussian random vectors are correlated.
Without loss of generality, both $X_A,Y_A$ can have zero-mean distributions. 
The covariances of the Gaussian variables $C_X(\alpha,\alpha')$ and $C_Y(\alpha,\alpha')$ are assumed to be functions only of $\alpha - \alpha'$. 
The choice of a truncation map parametrization is made as a trade-off between simplicity related to the computational cost and efficiency related to the model estimation given some data.
\begin{figure}[!ht]
\begin{tabular}{ccc}
\includegraphics[width= 0.3\textwidth,trim= 65 0 85 0, clip=true]{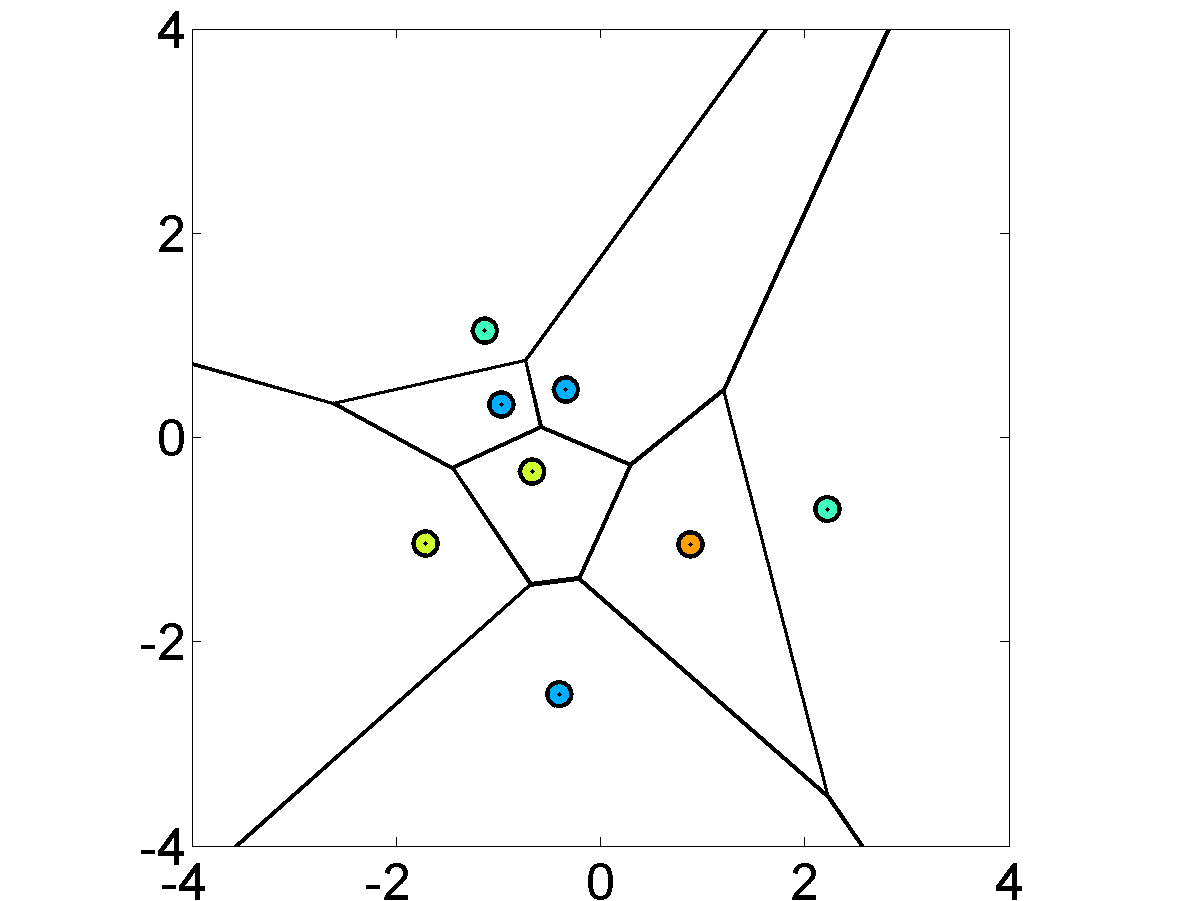}
&\includegraphics[width= 0.3\textwidth,trim= 65 0 85 0, clip=true]{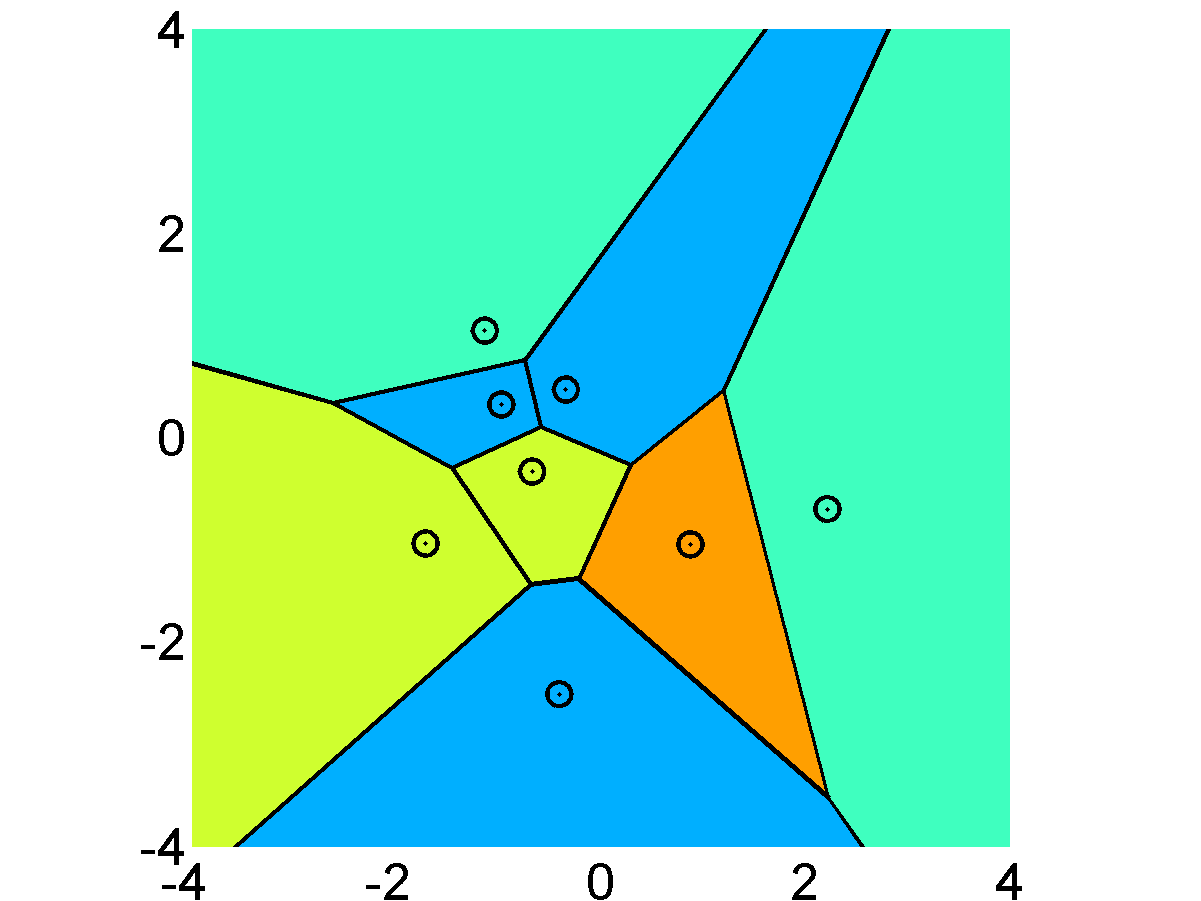}
&\includegraphics[width= 0.3\textwidth,trim= 65 0 85 0, clip=true]{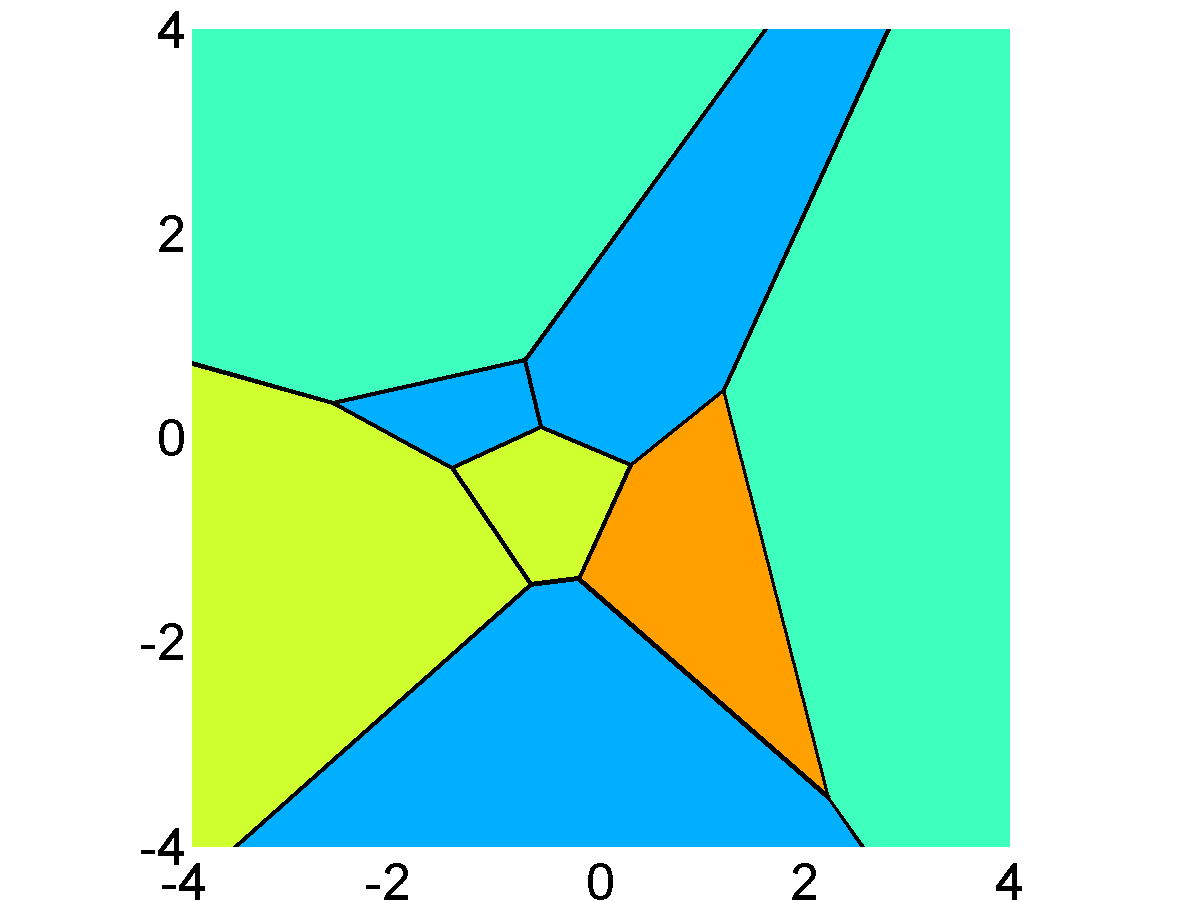}
\\
\parbox[b][2em][t]{0.3\textwidth}{
\footnotesize (a) Nodes with their categories} & 
\parbox[b][2em][t]{0.3\textwidth}{
\footnotesize (b) Coloring of $\mathbb{R}^2$ with category of the closest node} &
\parbox[b][2em][t]{0.3\textwidth}{
\footnotesize (c) Truncation map representation in $[-4,4]^2$} \\
\end{tabular}
\caption{The representation and interpretation of a truncation map based on colored Voronoi tessellation}
\label{fig:ex_tr_map}
\end{figure}
The representation and the interpretation of the truncation map is shown on the example Fig.~\ref{fig:ex_tr_map} of a sample from a prior distribution with mean number of nodes equal to 8. Figure~\ref{fig:ex_tr_map}(a) gives the positions and the colors (categories) related to the nodes, together with the borders of the Voronoi cells. Then the categories are assigned to every value in $\mathbb{R}^2$, which results in coloring Fig.~\ref{fig:ex_tr_map}(b). As far as every latent pair has standard Gaussian distribution, only the practical domain $[-4,4]^2$ is presented. The nodes are omitted in Fig.~\ref{fig:ex_tr_map}(c), while the interpretation of the mapping of a latent pair $(x,y)$ would be the color of the point with the coordinates $(x,y)$.
If greater simplicity of the partitioning is desired, other tessellations that reduce geometrical complexity might be considered.  
For instance, Voronoi tessellation with the $L_1$ metric limits the inclination of the tessellation borders to diagonal, vertical and horizontal. 
Another possibility is the T-tessellation \citep{kieu:13} which can be specified to have rectangular categorical areas, or even allow horizontal and vertical borders only. 
Even a simple coloring of cells with fixed thresholds can be considered. 
In terms of model estimation, all these types of parametrization could be used within the approach of this work, requiring only that the set of the prior distributions of the tessellation parameters are modified appropriately.

\subsubsection{Prior information on data parameters}
The model above is a probabilistic model. When the covariance functions $C_X,C_Y$ are fixed, the uncertainty in the model is due to the prior distribution on the parameters $\theta$ and due to the distribution of the latent random Gaussian vectors $X_A,Y_A$. 
When, in addition, $\theta$ are fixed, all categorical realizations $M_{\theta}(X_A,Y_A)$ follow the same distribution, 
and the model randomness is a result only of different realizations $x_A,y_A$.
When using stationary GRFs to generate categorical variables, the categorical field from the mapping function is also stationary.
Because of the stationarity, the distribution of several 
$z_B, B \subset A$ 
might be experimentally available. 
Denote as $B_h$ some bivariate vector $(\alpha_1,\alpha_2) \subset A$, such that $\alpha_1 - \alpha_2 = h$.
In this work, a particular case is considered. Distribution of $z_{B_h}$ is assumed to be known for $h=(1,0)$ and $h =(0,1)$ and some conventional unit length. It will be called bivariate unit-lag distributions. 
This is rather realistic case where several samples for the empirical distribution of $z_{B_h}$ may be available along a vertical or a horizontal well, respectively. With no restriction, a three-dimensional case would have been treated similarly. The empirical probability of $z_{B_h}$ is denoted $\pi(z_{B_h})$.

\begin{figure}[!ht]
\begin{center}
\begin{overpic}[width=0.5\textwidth]{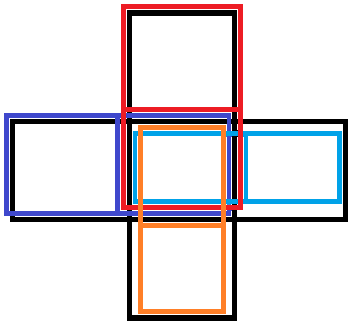}
\large
\put(11,64){\textcolor{blue}{$ \in {\cal B}$}}
\put(0,0){}
\put(71,75){\textcolor{red}{$\in {\cal B}$}}
\put(0,0){}
\put(74,18){$B_*$}
\put(0,0){}
\put(78,43){\textcolor{cyan}{$\in {\cal B}$}}
\put(0,0){}
\put(45,15){\textcolor[rgb]{1.0, 0.49, 0.0}{$\in {\cal B}$}}
\end{overpic}
\caption{A five-point patterns $B_*$ with known marginal distribtions on elements of ${\cal B}$} 
\label{fig:pattern}
\end{center}
\end{figure}
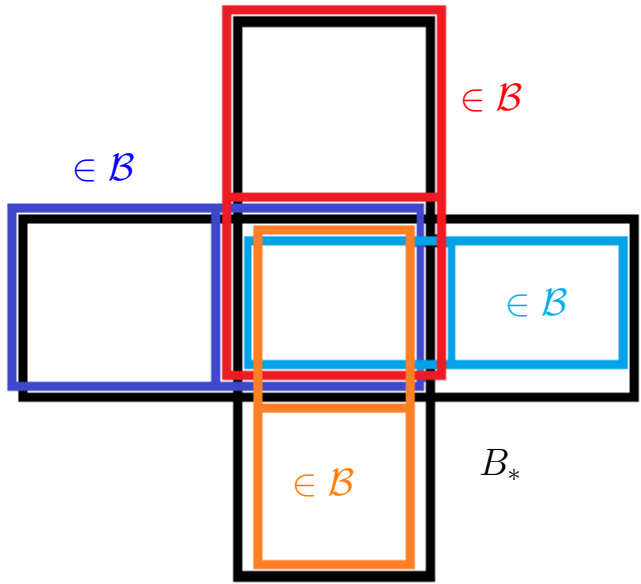

The truncation map parameter estimation should be based on this data. In order for the bivariate distributions to be consistent with each others and to create further on a likelihood function to account for both of them, the distributions will be
integrated in one joint distribution of a few observations (a pattern). The pattern distribution is then assumed correct in the truncation map estimation process.
While the joint pattern distribution retains maximum uncertainty,
the bivariate distributions represent marginal distributions for the pattern distribution. The pattern distribution is estimated even in case some of the bivariate probabilities are missing. The choice of pattern made in this work is a five-point neighborhood. It is denoted $B_* \in A$, $B_*=(\alpha_1,\dots,\alpha_5)$ is a sequence of length five, where the points form a five-point neighborhood, one point in the middle, and four surrounding points. The bivariate subsets of $B_*$ with known distributions will be denoted ${\cal B}$. The black frame in Fig.~\ref{fig:pattern} shows the pattern and the frames of other colors show the known marginal bivariate distributions on ${\cal B}$. The pairs of blue and light-blue elements, as well as the pairs of red and orange and elements have the same distribution, respectively due to stationarity.

Then the  problem of finding a probability mass function (pmf) $p(z_{B_*})$ that maximizes the entropy 
\begin{equation}\label{eq:entropy}
H(p) = - \sum_{z_{B_*}} p (z_{B_*}) \, \ln p (z_{B_*}) 
\end{equation}
is addressed.
This maximization problem is known as the Bayesian maximum-entropy (BME) principle. 
Here the distribution of the pattern $B_*$ is conditional to the marginal distributions on $\mathcal{B}$.
The entropy as given above is a measure of the amount of uncertainty represented by a probability distribution. \citet{jaynes:57} states that the BME estimate is the maximally noncommittal estimate with regard to missing information.

The existence of an optimal distribution on  $B_*$ is guaranteed as soon as the set of distributions that meet all 
constraints is nonempty. 
This is typically the case when all prespecified pmf derive from the 
same empirical distribution based on a training image on $A$ or well data including categorical observations at a given lag. 
Even though, the optimal distribution is not necessarily unique, one optimal pmf $p_*$ on $B_*$ can be found in order to further proceed with maximum likelihood estimation of parameters $\theta$ based on the data given by $p_*$.
The mathematical manipulations below are simplified using the notation $z=z_{B_*}$. 
To compute the distribution that maximizes the entropy (Eq.~\eqref{eq:entropy}) subject to 
the marginal distributions $\pi(z_B)$ for $B\in {\cal B}$, and the requirement that 
probabilities must sum to one, a Lagrangian is formed 
 for $p^*$, the Lagrangian is formed, 
\begin{multline}\label{eq:lgn}
L = - \sum_{z} p (z) \ln p (z) + \sum_{B \in {\cal B}} \sum_{z_B \in C^{|B|}} \lambda (z_B) \left[ \sum_{z_{B_*\backslash B}} 
p (z_B, z_{B_*\backslash B}) - \pi (z_B) \right] \\+ \lambda \left[ \sum_{z} p (z) - 1 \right]. 
\end{multline}
where the $\lambda (z_B)$'s and $\lambda$ are Lagrange multipliers. Maximizing the Lagrangian leads to the 
linear system of equations
\begin{equation}
 \left\{ 
\begin{array}{rll} 
\dfrac{\partial L}{\partial p (z)} & = & - \ln p (z) -1 + \sum_{ B \in {\cal B}} \lambda (z_B) + \lambda = 0; \\[0.4cm]
\dfrac{\partial L}{\partial \lambda (z_B)} & = &  \sum_{z_{B_*\backslash B}} 
p (z_B, z_{B_*\backslash B}) - \pi (z_B) = 0; \\[0.4cm]
\dfrac{\partial L}{\partial \lambda } & = & \sum_z p (z) - 1 = 0. 
\end{array}
\right. 
\end{equation}
The first set of equations gives that $p^*$ belongs to the exponential family
\begin{equation}
 p^* (z) = e^{\lambda -1} \prod_{B \in {\cal B}} e^{\lambda (z_B)}. 
\end{equation}
Moreover, using $1 = \sum p^* (z) $, we obtain
\begin{equation}\label{eq:pdev}
p^* (z) = \dfrac{ \prod_{B \in {\cal B}} e^{\lambda (z_B)} }{ \sum_{z} \prod_{B \in {\cal B}} e^{\lambda (z_B)} },
\end{equation}
as well as a formula relating $\lambda$ to all other Lagrange multipliers
\begin{equation}\label{eq:lambda}
\lambda = 1 - \ln \left( \sum_{z} \prod_{B \in {\cal B}} e^{\lambda (z_B)} \right).
\end{equation}
Solving this nonlinear system of equations is not easy, but an alternative method is available. 
By analogy with an algorithm proposed by \citet{deming:40} to estimate the entries of an array with 
prespecified marginal totals,  the following iterative algorithm to estimate $p^*$ is proposed. 

\begin{algorithm}

\caption{Deming and Stephan algorithm outline}\label{agm:der}

(i) initialize $p^c$; 

(ii) generate $B \sim {\cal U} ( {\cal B} )$; 

(iii) put $ p^n (z) = p^c (z) \pi (z_B) / p^c (z_B)$ for each $ z \in C^{|B_*|}$; 

(iv) put $p^c = p^n$ and goto (ii). 

\end{algorithm}

Its principle is quite simple. At each iteration, a family of indices $B \in {\cal B}$ is selected at random. 
Then the current distribution is updated to be assigned the prespecified $B$-margin. 
Unfortunately this latter step affects all other marginal distributions. The algorithm must be made iterative so that all $B$-marginals are respected in the long run. Eventually, this algorithm can be seen as a successive projections algorithm, one projection consisting of fixing a $B$-marginal. Such a method has an interesting geometric interpretation. To fix ideas, let ${\cal M}_B$ be the set of all distributions on $C^{|B_*|}$ with fixed $B$-marginal $\pi (\cdot_B)$. 
Let also $p$ be an arbitrary distribution on $C^{|B_*|}$. Then it can be shown that the distribution $p^*$ defined as 
\begin{equation}\label{eq:defpst}
p^* (z) = p (z) \dfrac{\pi (z_B)}{p (z_B)} 
\end{equation}
satisfies 
\begin{equation} \label{eq:min_KL}
D ( p^* \|  p ) = \min_{q \in {\cal M}_B} D ( q \| p ) , 
\end{equation}
where $ D ( q \| p )$ denotes the Kullback-Leibler divergence of $q$ w.r.t. $p$: 
\begin{equation} 
D ( q \| p ) = \sum_v q (z) \ln \dfrac{ q (z)}{p (z)} 
\end{equation}
This optimal distribution $p^*$ is called the {\it I-projection} of $p$ to ${\cal M}_B$. 
\citet{csiszar:75} proposed a criterion for characterizing $I$-projections. $p^*$ is called the $I$-projection 
of $p$ on ${\cal M}_B$ if and only if 
\begin{equation}\label{eq:czr}
D ( q \|  p ) \geq D ( q \|  p^* ) + D ( p^* \|  p ) \qquad  q \in {\cal M}_B. 
\end{equation}
One may check that $p^*$ as defined above by (\ref{eq:defpst}) satisfies (\ref{eq:czr}). Indeed, the following equality 
\begin{equation}\label{eq:decomp}
D ( q \|  p ) = \sum_{z_B} q (z_B) D \bigl( q ( \cdot_{B_*\backslash B} | \cdot_B ) \|  p ( \cdot_{B_*\backslash B} | \cdot_B ) \bigr) + 
D \bigl( q ( \cdot_B ) | p ( \cdot_B ) \bigr)
\end{equation}
always holds. Moreover, Eq.~\eqref{eq:czr} is satisfied by $p^*$ as an equality. One can readily see that 
\begin{align}
 \sum_{z_B} q (z_B) D \bigl( q ( \cdot_{B_*\backslash B} | \cdot_B ) \|  p ( \cdot_{B_*\backslash B} | \cdot_B ) \bigr) = D ( q \| p^* ); 
\\  D \bigl( q ( \cdot_B ) \| p ( \cdot_B ) \bigr) = D ( p^* \|  p ). 
\end{align}
A proof of the equalities is given in Appendix A. 
The estimated distribution on the pattern $B_*$, that integrates the known marginal distributions on $B \in {\cal B}$ is further used ot estimate the truncation map parameters. The BME estimate on $B_*$ is assumed to be the correct data in this following estimation.

\subsubsection{Truncation map goodness of fit}

The model parameters $\theta$  should be estimated based on their prior distribution and the data. 
The prior distribution of the parameters of a truncated bigaussian 
model based on Voronoi tessellation were given as distributions on the 
number of nodes, the locations of nodes, and the colors of nodes (see 
discussion in Sect.\ref{sec:pr_info_mod_par}).
The search for parameters of the truncation map should result in realizations of truncated bigaussians with approximately the same probability of occurrence of $z_{B_*}$ as given by the BME approach. 
Notice, that the correlations of the latent random Gaussian vectors on $B_*$ are assumed to be known. Samples with zero-mean and known correlation matrices of the latent Gaussian realizations on $B_*$ are denoted $x_{B_*}^{(1)}, \cdots ,x_{B_*}^{(n)}$ and $y_{B_*}^{(1)}, \cdots ,y_{B_*}^{(n)}$, respectively for some $n \in \mathbb{N}$.
The categorical relative frequency from the mapped sample is denoted 
\begin{equation}
\hat{f}(z_{B_*},\theta) = \frac{1}{n}\bigl|\{ M_{\theta} (x_{B_*}^{(r)},y_{B_*}^{(r)}) = z_{B_*},
 \quad r=1,\dots,n \} \bigr| 
\end{equation}
where $|\cdot|$ denotes the cardinality of the set.

The mismatch between the pmf $p^*$ and the observed frequency can be quantified using a divergence function such as the Kullback-Leibler divergence 
\begin{equation}
F(\theta) = \sum_{z_{B_*}} \hat{f}(z_{B_*},\theta)  \ln \frac{\hat{f}(z_{B_*},\theta)}{p^*(z_{B_*})}.
\label{eq:likel_f}
\end{equation}
The Kullback-Leibler divergence  is non-negative, taking the value zero only when the distribution of samples is equal to the expected distribution.  The Kullback-Leibler divergence is seen above to have a close relation to the BME estimates Eq.~\eqref{eq:min_KL}, although other divergence functions \citep{bregman:67} with similar properties could be used instead. It is assumed that whenever a facies included in the pattern realizations $z_{B_*}$ is not presented in the truncation map presented by the model parameters $\theta$, the mismatch is infinitely large. This consequently gives the estimates with each facies assigned to at least one cell almost surely.

\subsubsection{Simulated annealing algorithm}

The problem is set as maximum likelihood estimation represented by the mismatch function $F(\theta)$ in \eqref{eq:likel_f} above with negative sign. 
Parameters that minimize the mismatch function are estimated
 by means of the simulated annealing algorithm. The target distribution at temperature $T$ can be written as
\[
\mathcal{L}(\theta | p^*) = C(T) \exp(-\frac{F(\theta)}{T}),
\]
where temperature parameter $T$ changes with iterations, and $C(T)$ is an unknown normalization constant, dependent on $T$.
The temperature cooling schedule $T^{(n)}$, where $n$ is the iteration number is chosen to ensure that the algorithm does not get stuck in a local minimum at the early iterations. The temperature is typically allowed to cool slowly, for example, $T^{(n)} = T^{(0)} \alpha^n$ for $\alpha < 1$.
The earlier period of the iterative process have the larger temperature values $T$ to better explore the $\theta$ space.
The cooling period of the iterative process corresponding to the smaller values of $T$ would reach a local minimum of the mismatch. 
This corresponds to a local maximum of the goodness of fit of  the truncation map $M_{\theta}$ to the the data expressed as the BME estimate. 
The candidate-generating distribution in the iterative process is based on  the prior model in order to get an estimate close to the prior distribution. 
However, because the problem is stated as maximum likelihood estimation, the estimate does not necessary respect the distribution of the number of nodes $P_\mu(\cdot)$, the node colors or the node coordinates.
A standard implementation of the simulated annealing for maximization  of the function $F(\theta)$, with the the associated temperature cooling schedule $T(\cdot)$, and the candidate-generating distribution $q(\theta,\theta')$ is provided in Algorithm~\ref{agm:sim_a}.

\begin{algorithm}

\caption{Simulated annealing algorithm, given $\theta^{(0)}$, $T^{(0)}$, $q(\cdot,\cdot)$, $T^{(n)} = T^{(0)} \alpha^n$, $F(\cdot)$}
\label{agm:sim_a}

(i) initialize $\theta_c = \theta^{(0)}$; $T_c=T^{(0)}$

(ii) generate $\theta_n \sim q(\theta_c, \theta_n)$.

(iii) assign $\theta_c = \theta_n$ with probability 
\begin{equation}
	\rho (\theta_c,\theta_n;T_c) = 
	\min \left\{ 
	\frac{\exp (-F(\theta_n)/T_c)}{\exp (-F(\theta_c)/T_c)} , 1 \right\}.
\end{equation}

(iv) put $T_c = \alpha T_c$  and goto (ii).

\end{algorithm}

If $F$ is the objective (target) density and the acceptance probability $\rho(\theta_1,\theta_2, T)$ is replaced by 
\begin{equation}
	\rho (\theta_1,\theta_2) = 
	\min \left\{ \frac{\exp (-F(\theta_2))}{\exp (-F(\theta_1))} 
	\frac{q(\theta_1|\theta_2)}{q(\theta_2|\theta_1)}, 1 \right\},
\end{equation} 
which comes from setting $T=1$ and accounting for the probability of proposing a move,
this gives Metropolis-Hastings algorithm for sampling the likelihood.

The candidate-generating distribution $q$ based on the prior information on the model parameters is created as following. The truncation map is parametrized by a colored Voronoi tessellation. The prior distribution is conveniently specified in terms of a probability distribution of the number of nodes, distribution of the coordinates and the colors of the nodes. 
The proposals  for the  Markov chain includes  birth and death events. The proposed number of nodes $\nu_n$, given the current number of nodes $\nu_c$, is sampled from $\{\nu_c-1,\nu_c,\nu_c+1 \}$, with probabilities proportional to  $P_\mu(\nu_c-1),P_\mu(\nu_c)$ and $P_\mu(\nu_c+1)$ respectively,
where $P_\mu(\cdot)$ is the Poisson density with parameter $\mu$.
Depending on the outcome of the previous step, the proposal truncation map is sampled as following. 
In a birth event $(\nu_n = \nu_c +1)$ the location of  a new node is sampled at random according to the prior distribution and added to the current nodes. 
When the number of nodes stays the same $(\nu_c = \nu_c)$, an existing node is chosen uniformly and resampled according to the prior distribution.
In a death event $(\nu_n = \nu_c -1)$, one of the  existing nodes  is selected uniformly for removal.
Any of the events above produces a new tessellation, which might significantly change the map if the total number of nodes is small.

\subsection{Methods: Validation}

Any estimate of the truncation map together with the probability distribution of latent Gaussian random vectors provide a predictive distribution of categorical realizations at every location. 
Whenever observations at locations are available, this will be called an event. 
A good probabilistic model should provide adequate probability of generating the correct categories at those locations. 
The probabilities are only indirectly assimilated as marginals of the the pattern pmf and the bivariate pmf.
One one hand, the observations can be poor and spatially scattered to give the correct empirical estimates of the transition probabilities.  
One the other hand, the categorical observations may also include joint observations with longer correlations, which are not reflected in the data.
Moreover, the categorical observations can become available at a later time. 
Altogether, the validation method, described below can be applied due to any of the mentioned reasons as well as to evaluate any assumption of the estimation method. 
They include the following: the prior distribution of the truncation map parameters related to the Voronoi tessellation and comparison to different types of parametrization; hypotheses on the latent variable probability distribution, including covariance matrices; choice of the cooling schedule or the form of the goodness-of-fit function.
In this work, the quality of a probabilistic prediction will be measured using scoring rules. 
A scoring rule is a measure of the quality of a probabilistic prediction. 
It is said to be proper if the expectation of the scoring rule is maximized when the observations are drawn from the prediction distribution. 
The scoring rule of the TPG model parameters will be evaluated based on an ensemble of category predictions that materialize from a probabilistic model. 
In general, the comparison can only be made with observations or samples from the true distribution. 
However, two predictive distributions can be compared given the materialized event.

For an event that constitutes of a unique observation at one location, a variety of scoring rules are available \citep{gneiting:07a}. 
The sample space of this event is equal to the set of categories $C$ for the model introduced previously. 
Let the probabilistic forecast  be given as a pmf $p(\cdot)$ on $C$.
Then, one example of a scoring rule of $p$ given an event $c\in C$ is a logarithmic score,
\begin{equation}
S(p,c)= \log p(c),
\end{equation}
Conversely, for a vector of observations one might use so-called scoring rules of the model parameters given the unordered data. 
For the sake of simple notation the model parameters, denoted $\theta$,  
includes both truncation map parameter values, Eq.\eqref{eq:model_set_par} and the statistics of the latent Gaussian random vectors distributions, which govern together the predictive distribution.
A materialized event is denoted as $z_D=(z_d,d\in D)$ and is an observation vector at the locations $D \subset A$.
An example of a related scoring rule based on the logarithmic score \citep[][Eq.~(55)]{gneiting:07a} for the model parameters $\theta$,  given this event, is 
\begin{equation}
\label{eq:unordered_sc}
S_{\theta} = \sum_{d\in D}^n E_{(D\backslash d)} [ \log P (Z_d=z_d | Z_{(D\backslash d)}=z_{(D\backslash d)}, \theta)].
\end{equation}
Here at every location $d$ a logarithmic score is taken for the observation $z_d$ conditioned to the observations $V_{(D\backslash d)}$ at locations  $(D\backslash d)$. The set of locations  $(D\backslash d)$ is a random subset of all observations $D$ excluding $d$,
$(D\backslash d) \sim \mathcal{U}(\mathcal{P}(D \backslash d))$.  
The expectation is taken over $(D\backslash d)$, and those values are summed up for all locations $d\in D$.
In practice, one should approximate the probabilistic forecast for each event $Z_d=z_d  | Z_{(D\backslash d)}=z_{(D\backslash d)}, \theta$ of a categorical random variable $Z_d$ at location $d$ given the observations $z_{(D\backslash d)}$ at $(D\backslash d)$, and the model parameters $\theta$. 
The general steps of the simulation of a categorical observation $z_d$ conditional to the set of categorical observations $z_{(D\backslash d)}$ are included in Algorithm~\ref{agm:z_d_cond_unordered}.
\begin{algorithm}
\caption{Simulation of $z_d$ conditional to $Z_{(D\backslash d)}=z_{(D\backslash d)}$}\label{agm:z_d_cond_unordered}

(i)  Generate 
\begin{equation}
\label{eq:cond_sim}
(x_{(D\backslash d)},y_{(D\backslash d)})
\sim
{\cal N}\bigl( \bar{0}, 
\left[\begin{array}{cc}
C_{X_{(D\backslash d)}} & 0       \\
0 & C_{Y_{(D\backslash d)}}       \\
\end{array}\right]  | \theta, M_\theta (x_{(D\backslash d)},y_{(D\backslash d)}) = z_{(D\backslash d)} \bigr),
\end{equation} 
where $\bar{0}$ is a zero-vector of proper length, $C_{X_{(D\backslash d)}}$ and $C_{Y_{(D\backslash d)}}$ are the covariance matrices of $X_{(D\backslash d)}$ and $Y_{(D\backslash d)}$, respectively. Simulation of $(x_{(D\backslash d)}$ and $y_{(D\backslash d)}$ should be done simultaneously due to the conditioning.
This step represents itself an iterative procedure, such as the constrained Gibbs sampler or its alternative form described in Sect.~\ref{subsec:cond_sim} below.

(ii) Generate Gaussian observation $(x_d,y_d)$ at the location $d$ conditional to the vector of Gaussian observations at the locations $(D\backslash d)$, that is, simulate
\begin{equation}
x_d \sim {\cal N}\bigl( 0,1 | X_{(D\backslash d)}=x_{(D\backslash d)} \bigr),
\quad 
y_d \sim {\cal N}\bigl( 0,1 | Y_{(D\backslash d)}=y_{(D\backslash d)} \bigr),
\end{equation} 
where simulation of $x_d$, and $y_d$ can be done separately due to the independence assumption. 
A standard method for conditional Gaussian simulation is to use kriging estimates for the mean and variance.

(iii) Produce a categorical observation at $d$ using the mapping $M_\theta(x_d,y_d)=z_d$.
\end{algorithm}
The problem of conditioning random categorical fields produced from truncated Gaussian or plurigaussian random vector to several categorical observations is often necessary beyond the application to the validation with the scoring rules. While the second and the third steps are straightforward, the next section is dedicated to the first step, and the sampling method used particularly for the truncation map in form of a colored Voronoi tessellation.

\subsection{Methods: Conditional simulation}
\label{subsec:cond_sim}
A common practical tool for unconditional Gaussian simulation is the Gibbs sampler. 
The Gibbs sampler has also been used to simulate latent variables conditioned to categorical observations and governed by a bigaussian truncated model, despite the difficulty to compute or find the rate of convergence experimentally \citep{armstrong:11}.
\citet{lantuejoul:12} proposed another iterative version of the unconditional simulation, known as the propagative version of the Gibbs sampler. This sampler avoids  problems with the use of the moving search method and reduces the need for covariance matrix inversion.
The main idea is to operate the Gibbs sampler with $\tilde{x_A}=C_{X_A}\inv \, x_A$ instead of the direct updating of the realization $x_A$ of a random vector with zero-mean and the covariance matrix $C_{X_A}$. 
The initial vector $x_A^{(0)}$ in the approach has all the same entries (e.g., zero-vector).
Then at every update one index $\alpha \in A$ is chosen at random in such a way that all the indices are selected infinitely often as the number of iterations tends to infinity. 
For a given index $\beta \in A$ (pivot), an update  of $x_A^{(k-1)}$ at state $k-1$ to state $x_A^{(k)}$ for every component takes the form
\begin{equation} 
	x_\alpha^{(k)} = x_\alpha^{(k-1)} + (-x_\beta^{(k-1)} + u^{(k)})C_{X_A}(\alpha,\beta), \quad \forall \alpha \in A,
	\label{eq:alter_gibbs}
\end{equation} 
where $C_{X_A}(\alpha,\beta)$ is an element $(\alpha,\beta)$ of the covariance  matrix $C_{X_A}$ of $X_A$; $u^{(k)}$ is a standard Gaussian random realization independent of $X_A$.
An adaptation to the case of the conditional TPG simulation based on the propagative version of the Gibbs sampler was given by \citet{emery:14}. The author used the classical truncation map parametrization \citep{armstrong:11}, where the
problem of conditioning to categorical observations was equivalent to  the problem of simulation of a GRF realization $x_A$ with correlation matrix $C_{X_A}$, subject to linear inequality constraints 
\begin{equation} \label{eq:lin_eq_prob}
\underline{I}_\alpha \leq x_\alpha \leq \overline{I}_\alpha, 
\quad \underline{I}_\alpha \in [- \infty,\infty),\; \overline{I}_\alpha \in (- \infty, \infty], \quad \alpha \in A.
\end{equation} 
The relation of the linear inequality conditioning  on $x$ was translated to conditioning on $u$, and the initial state $x^{(0)}$ was chosen to be conditional to the categorical data.
This work extends the idea for the case of conditional Gaussian simulation based on truncation map with the new parametrization in form of colored Voronoi tessellation representing the truncation map. 
The main usefulness of the method consists in the improved mixing quality.

\begin{figure}[!ht]
\centering
\begin{tabular}{cc}
\includegraphics[width= 0.3\textwidth,trim= 50 0 65 0, clip=true]{1_map_ColVor} &
\includegraphics[width= 0.27\textwidth,trim= 50 0 65 0, clip=true]{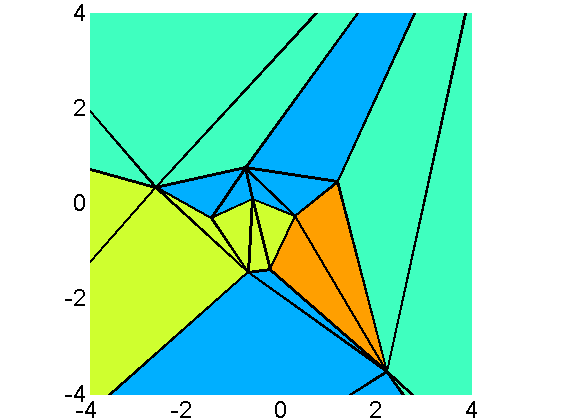} \\
\footnotesize (a) Truncation map & \footnotesize (b) Triangulation\\
\end{tabular}
\caption{Example of triangulation of a truncation map based on colored Voronoi tessellation}
\label{fig:ex_tr_map_tri}
\end{figure}

Conditioning to a categorical observation $z_\alpha$ within the model can be seen as conditioning of a pair of Gaussian realizations $(x_\alpha,y_\alpha)$, where the two-dimensional point $(x_\alpha,y_\alpha)$ belongs to the intersection of a union of polygons.
Due to possible triangulation of each cell of the Voronoi tessellation, 
one can write this subset as a union of triangles 
$\triangle_{t} \subset \mathbb{R}^2, t\in T(z_\alpha)$ 
with pairwise disjoint interiors. An example of a triangulation of the truncation map used in Fig.~\ref{fig:ex_tr_map} is shown in Fig.~\ref{fig:ex_tr_map_tri}. In order to triangulate only bounded Voronoi cells, a few additional nodes with large absolute coordinate values are added. Thus, some triangle edges are cut by the domain bounded by the practical range of the bivariate Gaussian distribution. Adding external nodes does not change the tessellation inside this domain. The truncation map is reproduced in Fig.~\ref{fig:ex_tr_map_tri}(a) for visualization purpose.

The equations for GRFs $(x_A,y_A)$ of the propagative version of the Gibbs sampler for a given pivot $\beta$, take form similar to Eq.~\eqref{eq:alter_gibbs}. 
Written explicitly for every element $\alpha \in A$, this gives the following.
\begin{equation}  \begin{split}
	x_\alpha^{(k)} & = x_\alpha^{(k-1)} + (-x_\beta^{(k-1)} + u^{(k)})C_{X_A}(\alpha,\beta)\\
	y_\alpha^{(k)} & = y_\alpha^{(k-1)} + (-y_\beta^{(k-1)} + v^{(k)})C_{Y_A}(\alpha,\beta),
\end{split}
\label{eq:alter_gibbs_gen}
\end{equation} 
where $C_{X_A}(\alpha,\beta)$ and $C_{Y_A}(\alpha,\beta)$ are the elements $(\alpha,\beta)$ of the covariance  matrices $C_{X_A}$ and $C_{Y_A}$ respectively. 
$u^{(k)}$ and $v^{(k)}$ are standard Gaussian random variables, independent of $X_A$ and $Y_A$.
Given the current states $x^{(k-1)}$, $y^{(k-1)}$, the updates in Eq.~\eqref{eq:alter_gibbs_gen} are functions of $u^{(k)}$ and $v^{(k)}$.
\begin{equation}  
	(x_\alpha^{(k)},y_\alpha^{(k)}) = 
	\Phi_{\alpha,\beta}(u^{(k)},v^{(k)},x_\beta^{(k-1)},y_\beta^{(k-1)},C_{X_A},C_{Y_A}) 
	\equiv 
	\Phi_{\alpha,\beta}(u^{(k)},v^{(k)}).
\label{eq:affine_transf}
\end{equation} 
$\Phi_{\alpha,\beta}$ is an affine transformation. 
Thus to update all $(x_\alpha,y_\alpha) \in \underset{t\in T(z_\alpha)}{\bigcup} \triangle_{t}, \alpha \in A$ 
at a given state $k$ given the pivot $\beta$, one should simulate 
\begin{equation} 
\label{eq:sim_point_uv_cond}
(u^{(k)},v^{(k)}) \in \bigcap\limits_{\alpha\in A} \bigcup\limits_{t\in T(z_\alpha)} \Phi_{\alpha,\beta}^{-1}(\triangle_{t}),
\end{equation}  
in order to satisfy the constraints. 
The sampling domain for $(u^{(k)},v^{(k)})$ can also be represented as a union of triangles with zero-probability of mutual intersections. 
Every iteration of the conditional simulation is straightforward from bivariate Gaussian simulation in a triangle and the cumulative distribution function over a triangle. The later is computed based on \citet[][Eq.~26.3.23 and Fig. 26.7--26.10]{zelen:12}.
Then a sampling of the two variables within a triangle is performed using a few iterations of the Gibbs sampler. 
The alternative sampler performs best when starting with several iterations of standard Gibbs to get considerable improvement from the initial state sampled as uncorrelated observations conditional to categories, then alternating standard Gibbs sampler and the alternative Gibbs sampler every iteration. The first one is still used for fast improvement of correlations while the second one allows to better explore the probable states through simultaneous updates of the entire latent vectors.
It has been found that a good initialization similar to the zero-state of the unconditional simulation, is to assign the same pairs of values to the pair of elements of the state 
$(x_\alpha^{(0)},y_\alpha^{(0)}) \; \forall \alpha \in A$ 
corresponding to the same categorical observation $z_\alpha = z\in C$, 
that is, 
\begin{equation}
\label{eq:ini_altern_Gibbs}
(x_\alpha^{(0)},y_\alpha^{(0)}) = (\tilde{x},\tilde{y})(z) \quad \forall \alpha \in A \; 
| \; z_\alpha = z = M_\theta((\tilde{x},\tilde{y})(z)),
\end{equation}
where the pairs of values $(\tilde{x},\tilde{y})(z)$ for $z\in C$ are chosen arbitrarily.

One scan of the sampler includes the updates of $(x_A,y_A)$, given every index $\beta \in A$ as pivot, where the order of choosing 
$\beta \in A$ is random.
One iteration of the sampler in the application to this work included one scan of a standard Gibbs sampler, and one scan of the conditional propagative version, which together provided fast convergence. Algorithm~\ref{agm:cond_sim_tpg} shows the steps of the iterative procedure.

\begin{algorithm}
\caption{Conditional simulation of the truncated bivariate model}\label{agm:cond_sim_tpg}

(i) Set $k=0$;

(ii) initialize $(x_A,y_A)^{(0)}$ according to Eq.~\eqref{eq:ini_altern_Gibbs};

\hfill Comment: propagative version scan

(iii) select a random path of indices $\beta \in A$;

\hspace*{5mm} (iii.a) 
for a selected $\beta \in A$ simulate $(u^{(k)},v^{(k)})$ according to Eq.~\eqref{eq:sim_point_uv_cond};

\hspace*{5mm} (iii.b) 
update $(x_\alpha^{(k)},y_\alpha^{(k)}) \; \forall \alpha \in A$ according to Eq.~\eqref {eq:alter_gibbs_gen};

\hfill Comment: standard Gibbs sampler scan

(iv) select a random path of indices $\alpha \in A$;

\hspace*{5mm} (iv.a) for a selected $\alpha \in A$ update 
\begin{equation}
\label{eq:cond_sim_1}
(x_\alpha^{(k)},y_\alpha^{(k)})
\sim
{\cal N}\bigl( 
(0,0), 
\left(\begin{array}{cc}
1 & 0       \\
0 & 1       \\
\end{array}\right) \; | 
\; X_{A\backslash\alpha} = x_{A\backslash\alpha}^{(k)},Y_{A\backslash\alpha} = y_{A\backslash\alpha}^{(k)}, 
\; \theta, 
\; M_\theta (x_\alpha^{(k)},y_\alpha^{(k)})  = z_\alpha \bigr),
\end{equation} 

(v) set $k = k+1$ and goto (iii).

\end{algorithm}

\section{Application}

\begin{figure}[!ht]
\begin{tabular}{cccc}
\includegraphics[width= 0.22\textwidth,trim= 5 0 5 0, clip=true]{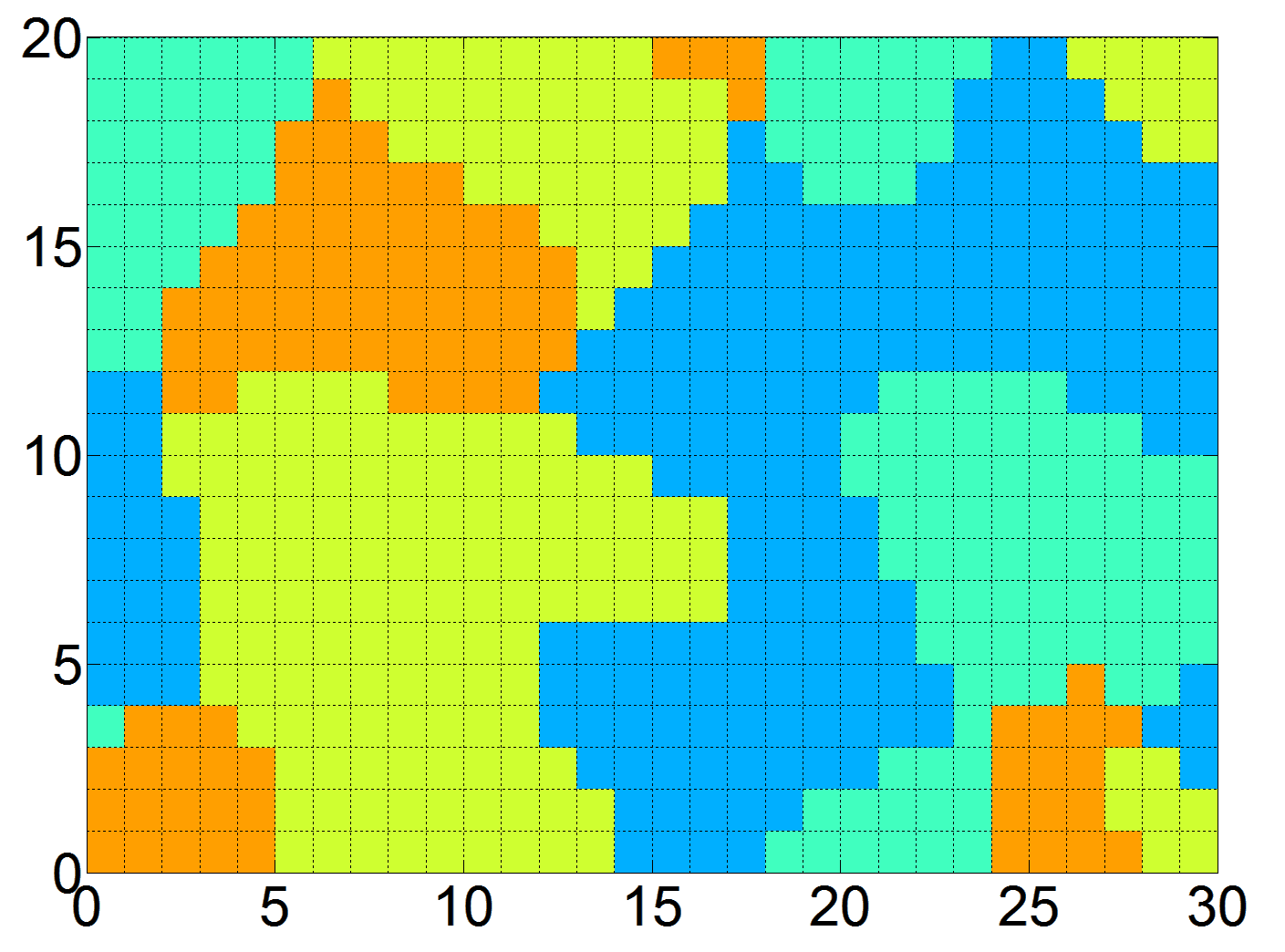} &
\includegraphics[width= 0.22\textwidth,trim= 5 0 5 0, clip=true]{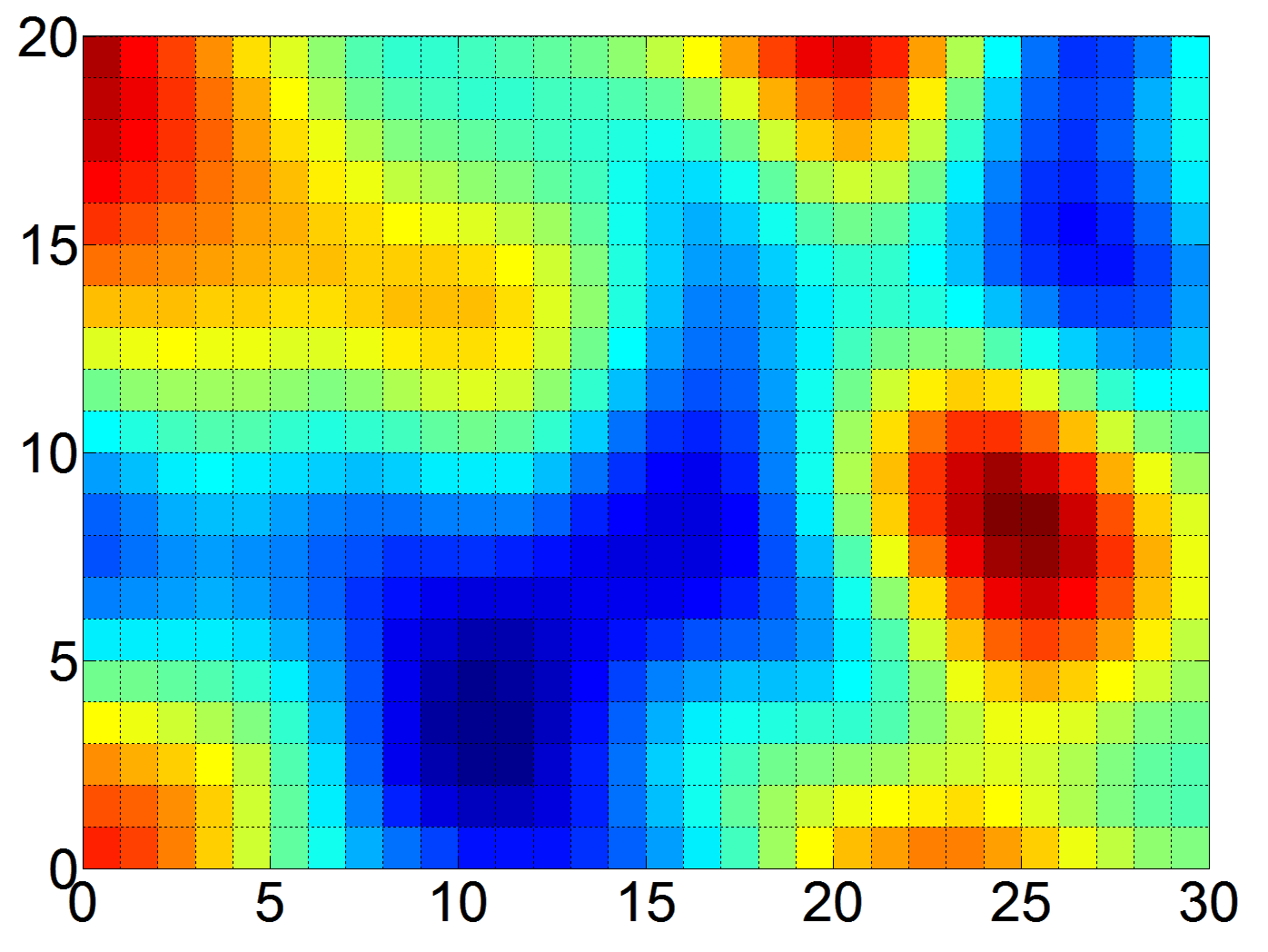}&
\includegraphics[width= 0.22\textwidth,trim= 5 0 5 0, clip=true]{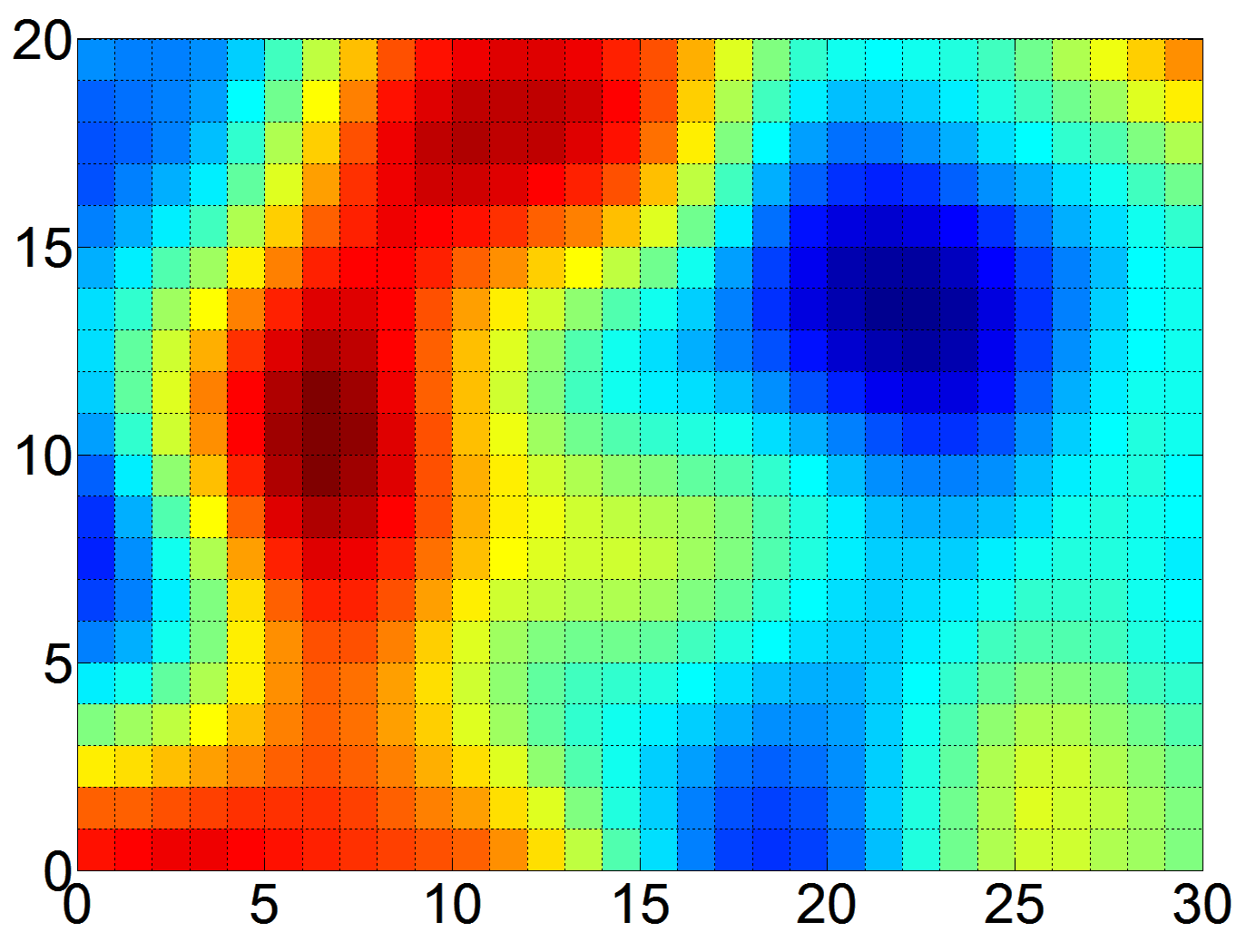}&
\includegraphics[width= 0.18\textwidth,trim= 45 0 65 0, clip=true]{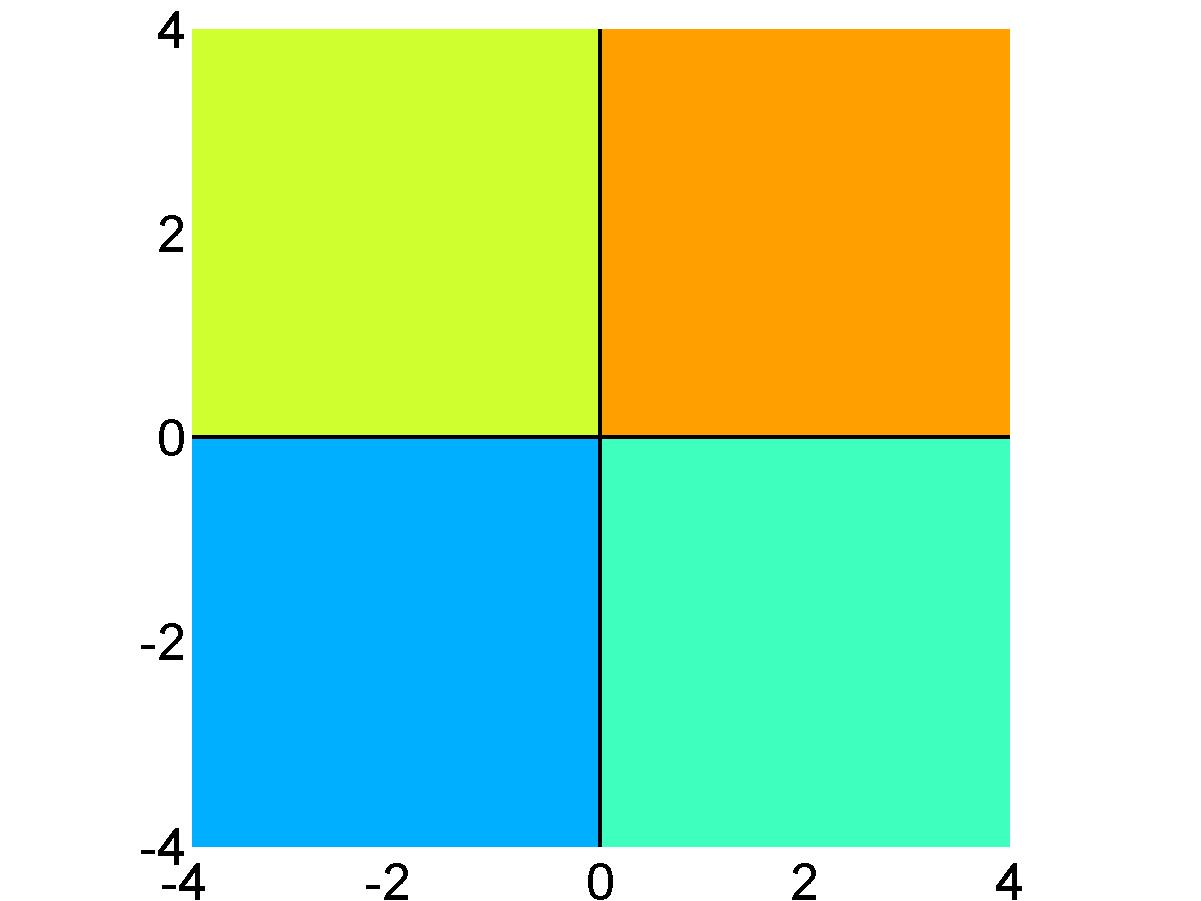} \\
\footnotesize (a) Synthetic field &
\footnotesize (b) Gaussian realization $x_A$ &
\footnotesize (c) Gaussian realization $y_A$ &
\footnotesize (d) Truncation map\\
\end{tabular}
\caption{The observations come from (a) which is created with the GRFs (b,c) and the truncation map (d)
\label{fig:synt_all}}
\end{figure}

The methodology for estimating a truncation map  from a realization of 
the categorical variables is illustrated with the following example.
A synthetic categorical field Fig.~\ref{fig:synt_all}(a) is generated from the 
truncation map in Fig.~\ref{fig:synt_all}(d) and from realizations of two latent Gaussian 
vectors (Fig.~\ref{fig:synt_all}(b,c)), each with a Gaussian covariance matrix with scale 
parameter equal to 10.

\subsection{Synthetic transition data}

\begin{figure}[!ht]
\begin{tabular}{cc}
\begin{overpic}[width= 0.3\textwidth,trim= 65 0 65 0, clip=true]{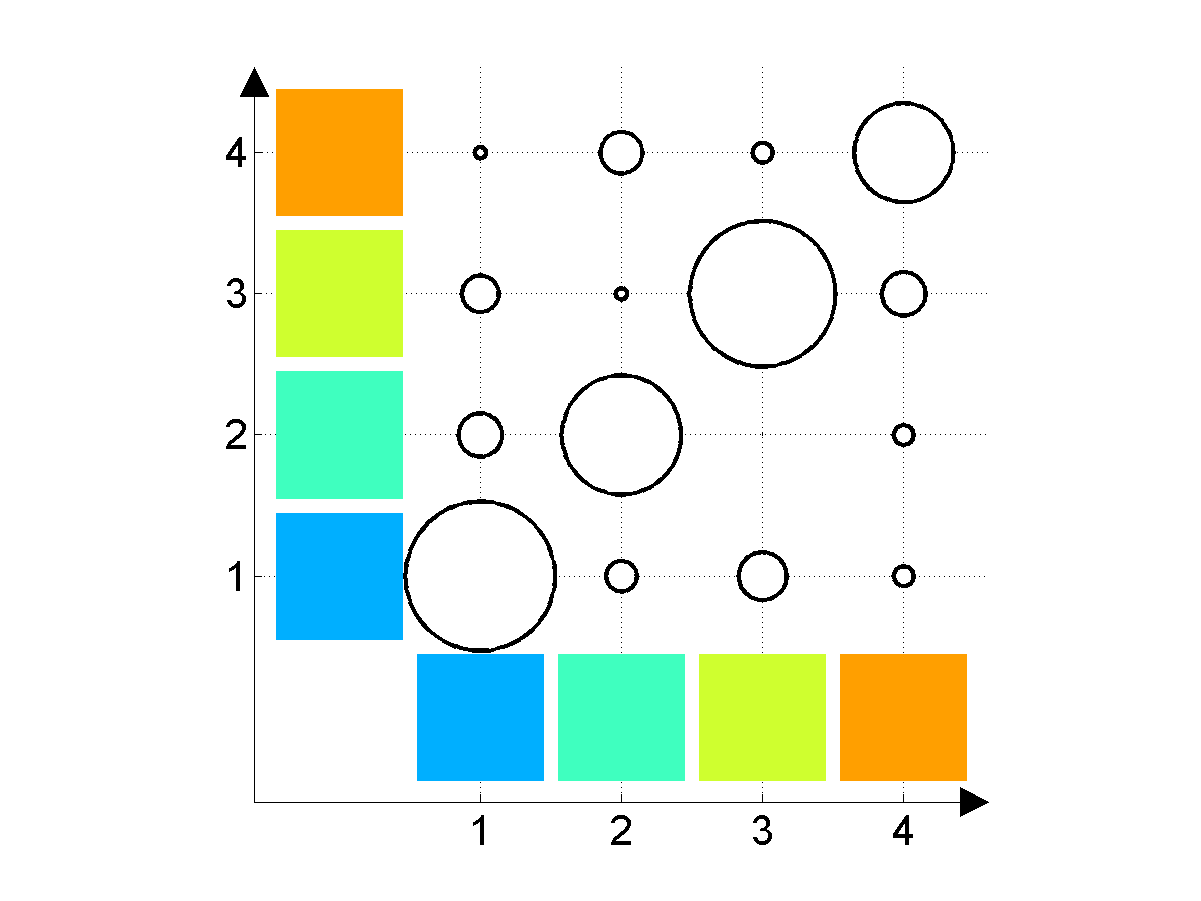}
\put(0,90){$z_{\alpha_2}$}
\put(0,0){}
\put(90,3){$z_{\alpha_1}$}
\end{overpic} 
&
\begin{overpic}[width= 0.3\textwidth,trim= 65 0 65 0, clip=true]{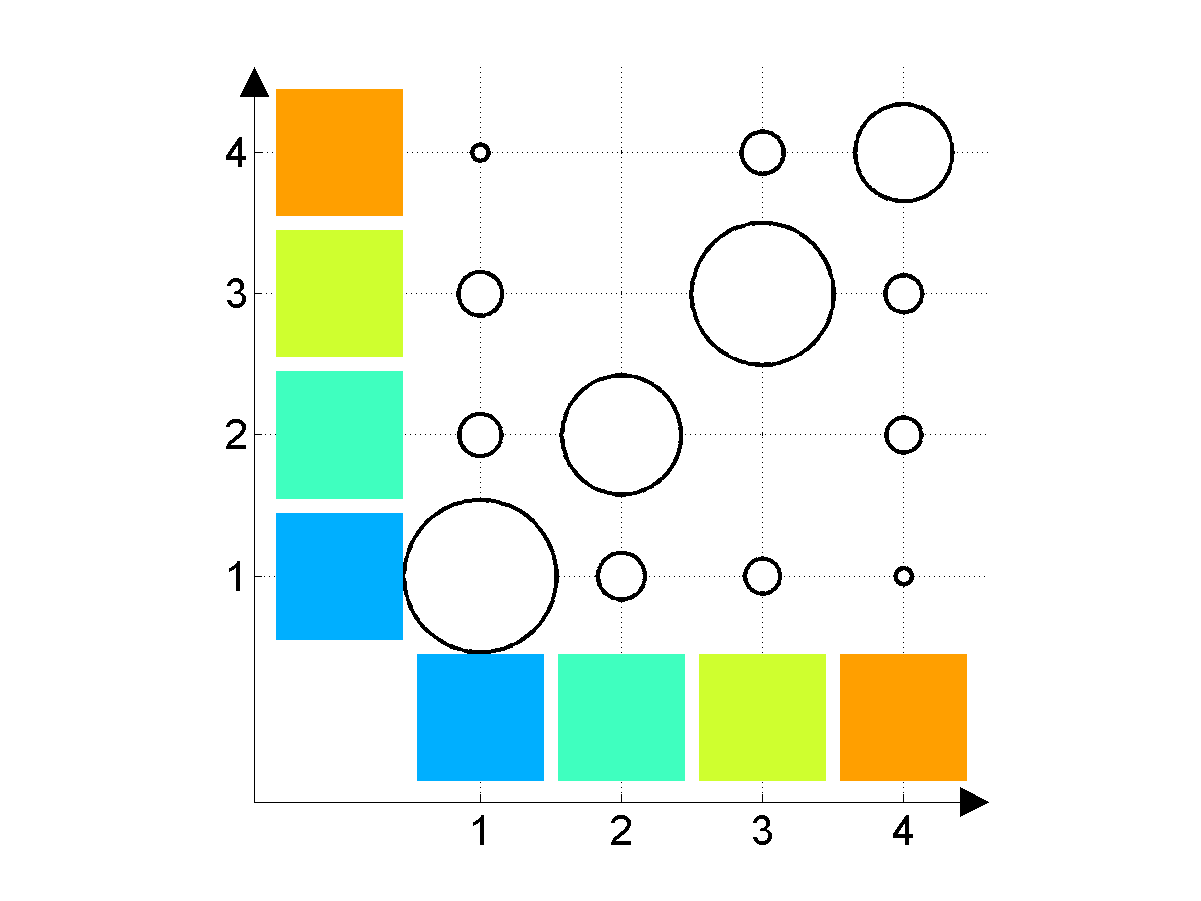}
\put(0,90){$z_{\alpha_2}$}
\put(0,0){}
\put(90,3){$z_{\alpha_1}$}
\end{overpic}\\
\footnotesize (a) 
$(z_{\alpha_1},z_{\alpha_2}) \in C^2$, 
$(\alpha_1,\alpha_2)=B \in \mathcal{B}_{(1,0)}$  
&
\footnotesize  (b) 
$(z_{\alpha_1},z_{\alpha_2})\in C^2$, 
$(\alpha_1,\alpha_2)=B \in \mathcal{B}_{(0,1)}$ \\
\end{tabular}
\caption{The area of circles at $(z_{\alpha_1},z_{\alpha_2})\in C^2, C=\{1,2,3,4\}$ is proportional to $\pi(z_{B})$ \label{fig:trans_data}}
\end{figure}
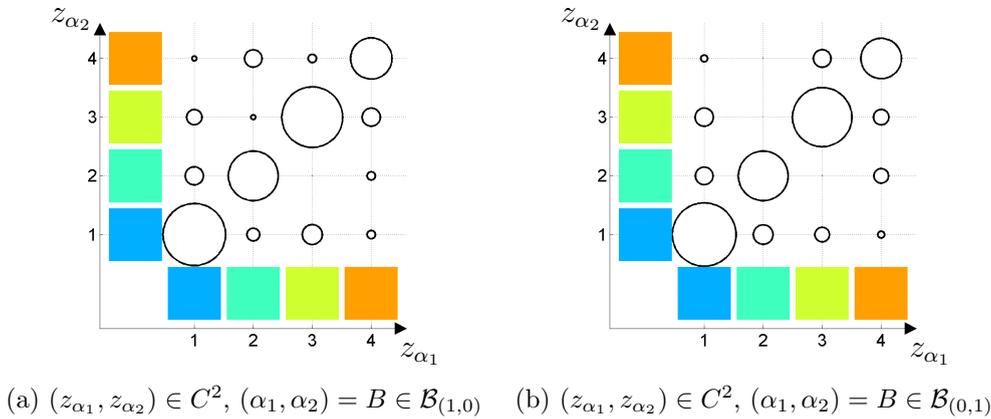

The transition probabilities that are used as data are derived from the synthetic categorical realization. 
For visualization, the transition probability for categorical realizations 
$(z_{\alpha_1},z_{\alpha_2}),\in C^2,C =\{1,2,3,4\}$, 
$(\alpha_1,\alpha_2)=B \in \mathcal{B}_h,h \in \mathcal{H}=\{(0,1), \ (1,0)\}$ 
are presented by areas of the circles with centers at the coordinates 
$(z_{\alpha_1},z_{\alpha_2})$ 
in Fig.~\ref{fig:trans_data}. 

\subsection{Estimation results}

\begin{figure}[!ht]
\begin{tabular}{rcccc}
 \raisebox{3ex}{\rotatebox{90}{First chain}} &
\includegraphics[width= 0.21\textwidth]{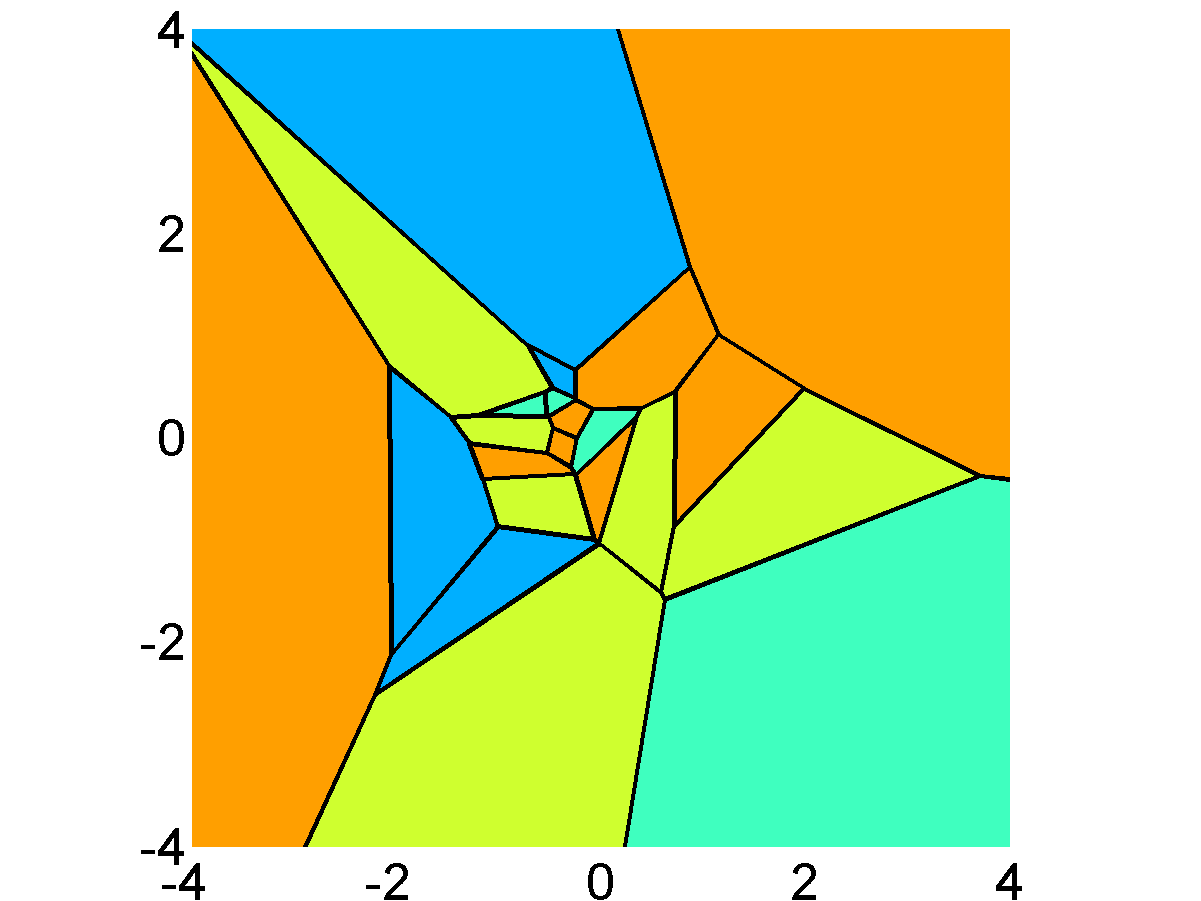} &
\includegraphics[width= 0.21\textwidth]{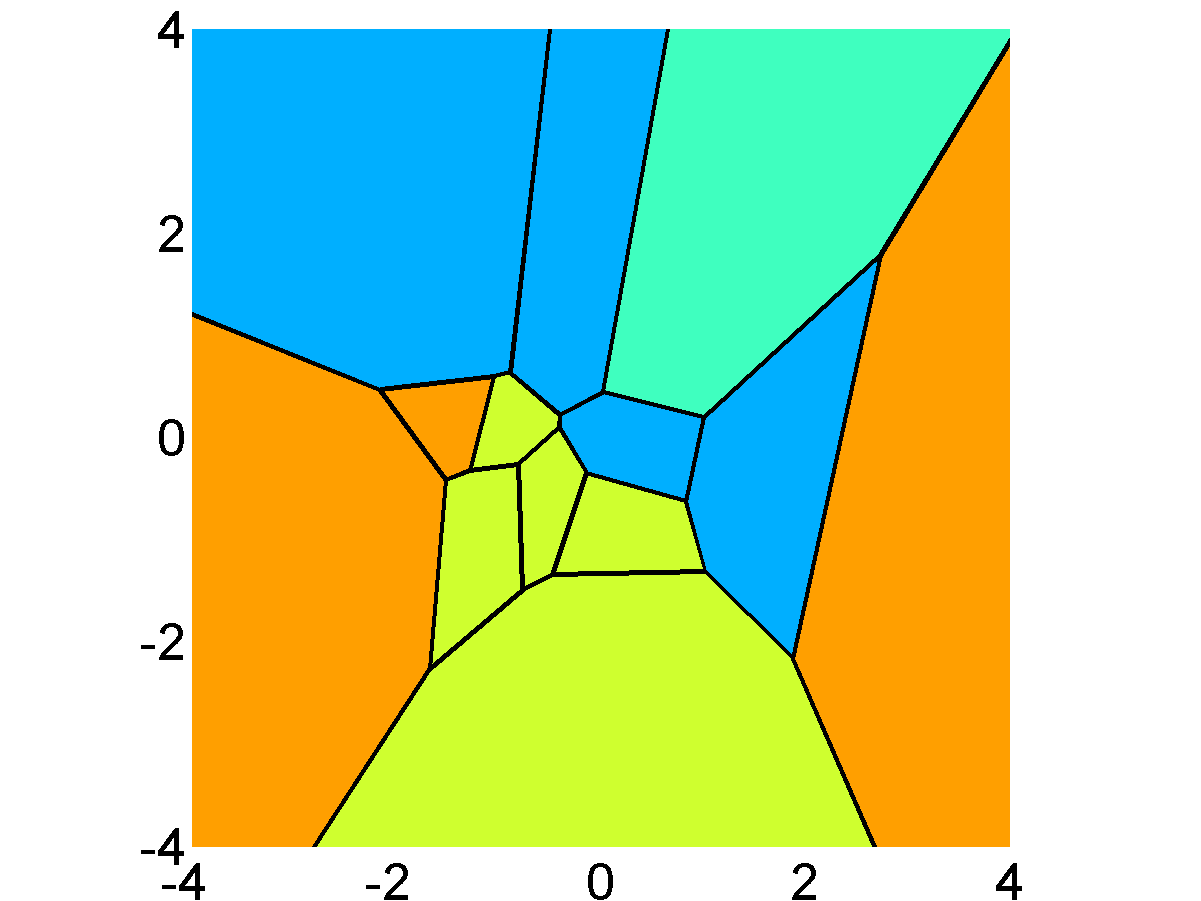} &
\includegraphics[width= 0.21\textwidth]{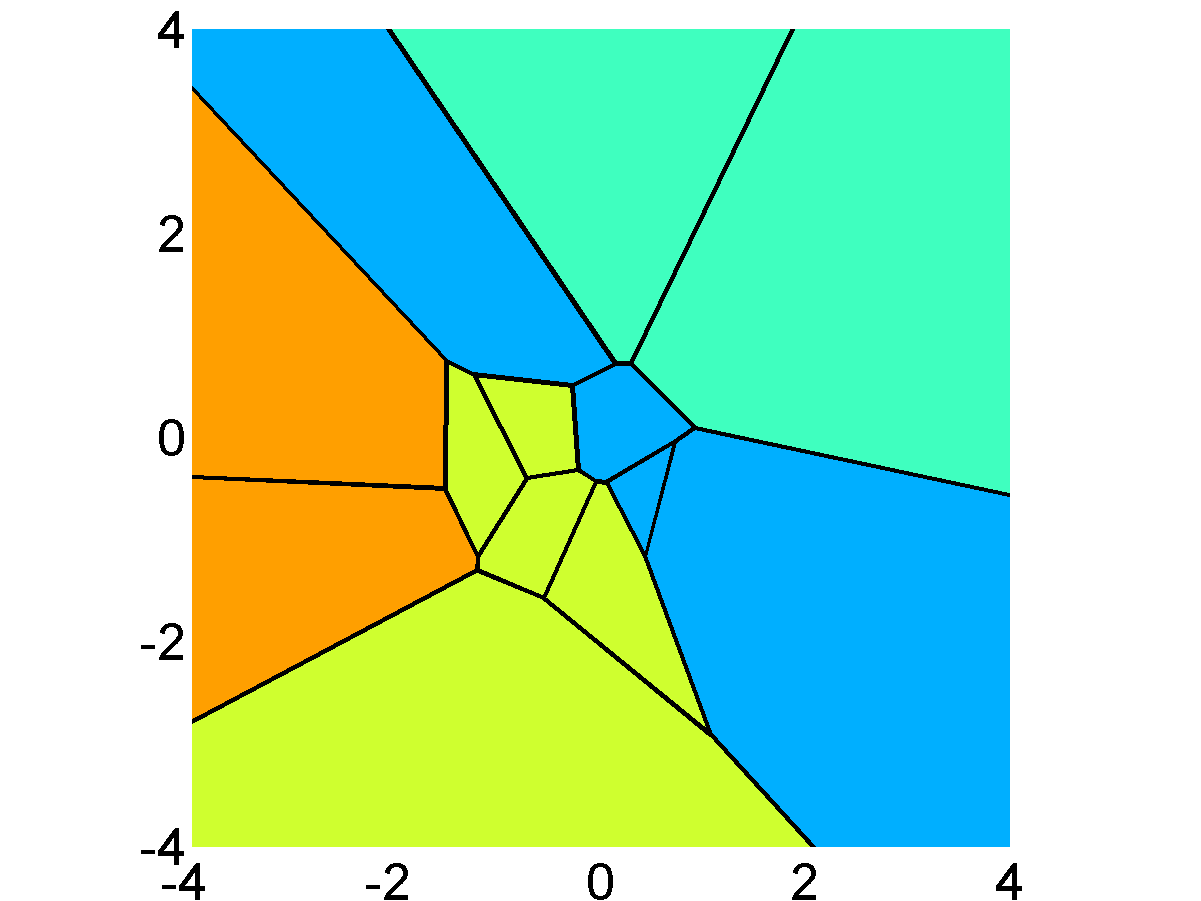} &
\includegraphics[width= 0.21\textwidth]{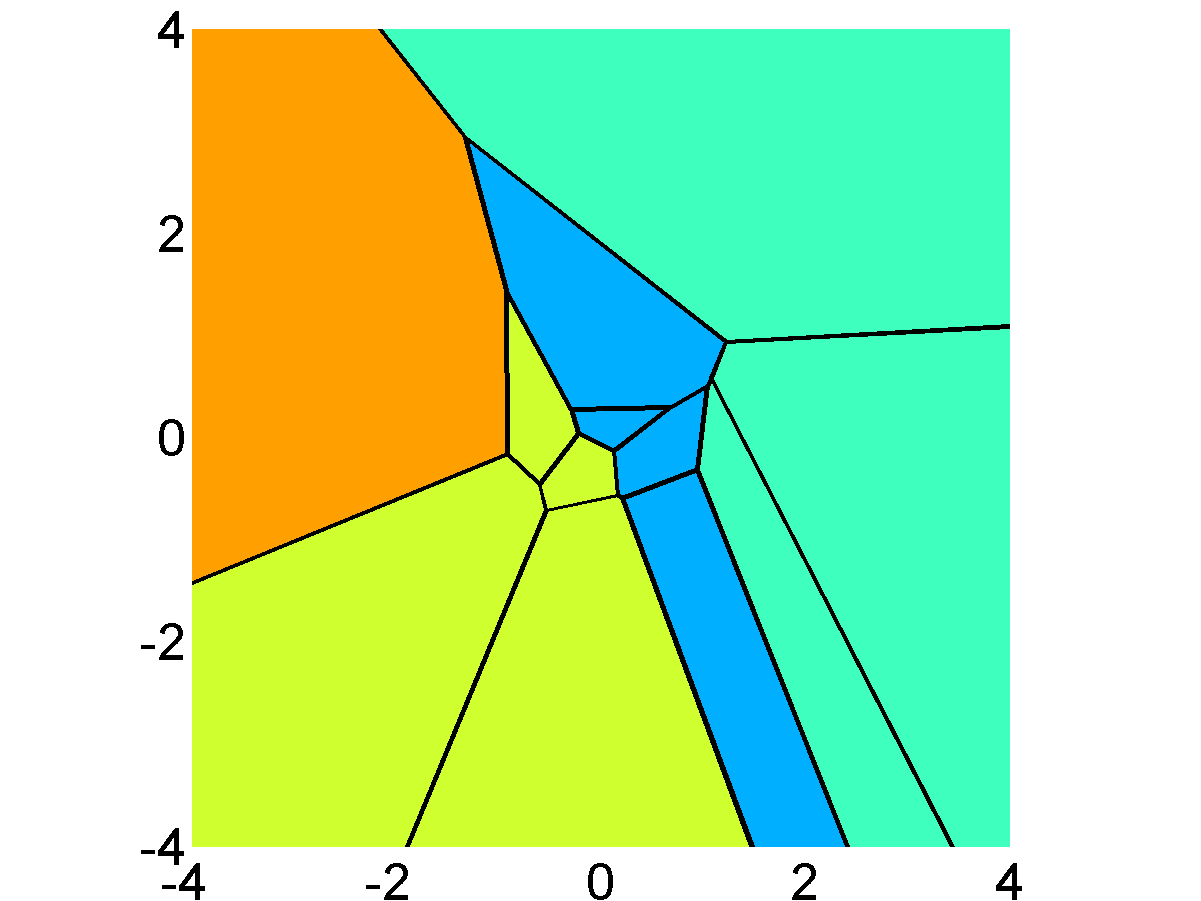} \\
 \raisebox{2ex}{\rotatebox{90}{Second chain}} &
\includegraphics[width= 0.21\textwidth]{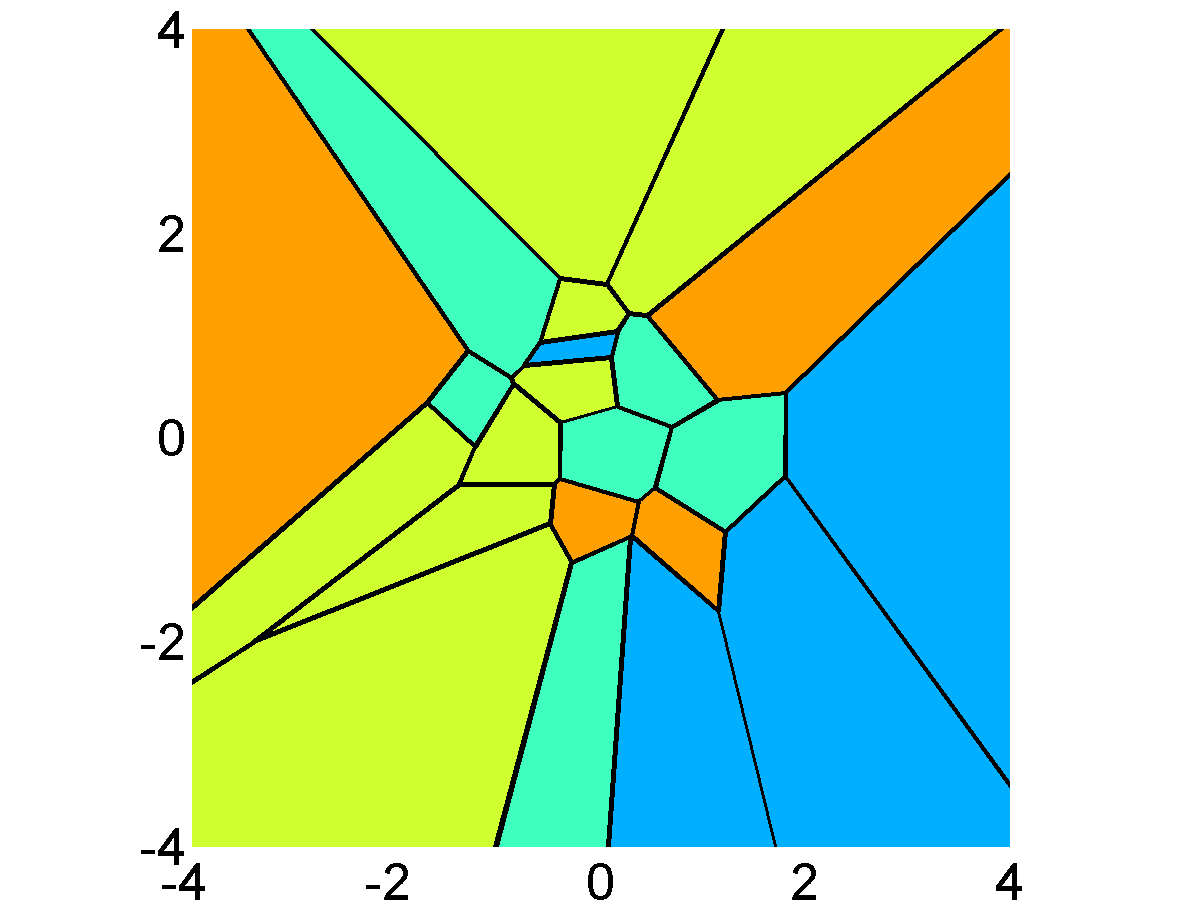} &
\includegraphics[width= 0.21\textwidth]{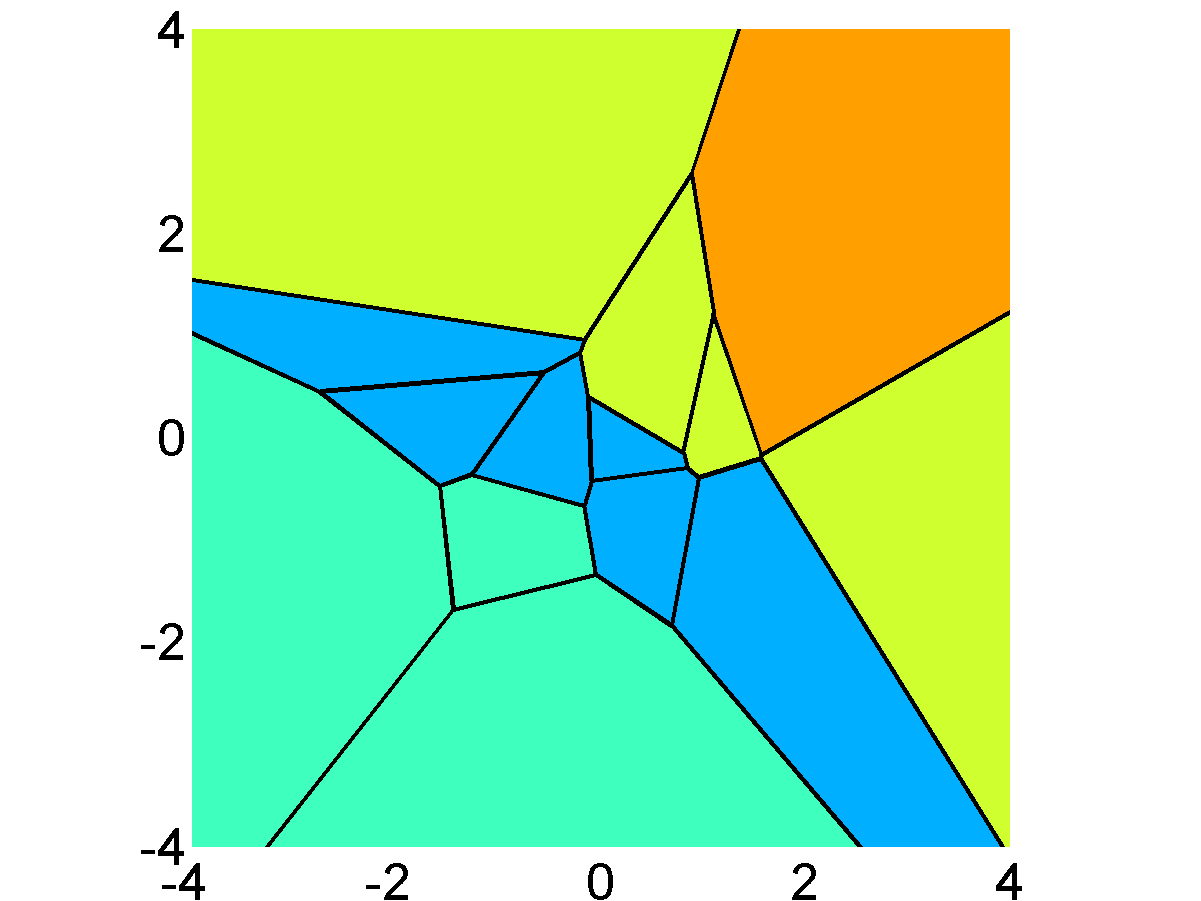} &
\includegraphics[width= 0.21\textwidth]{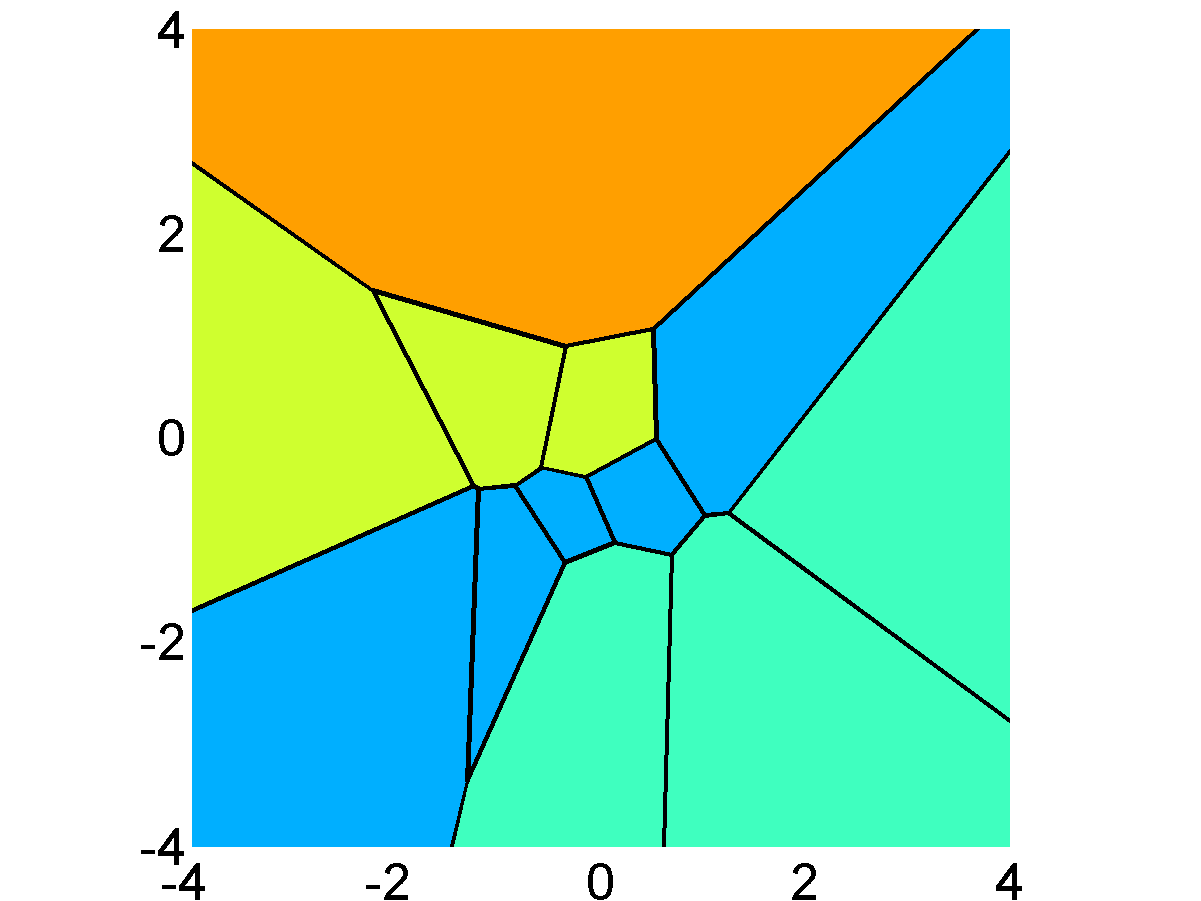} &
\includegraphics[width= 0.21\textwidth]{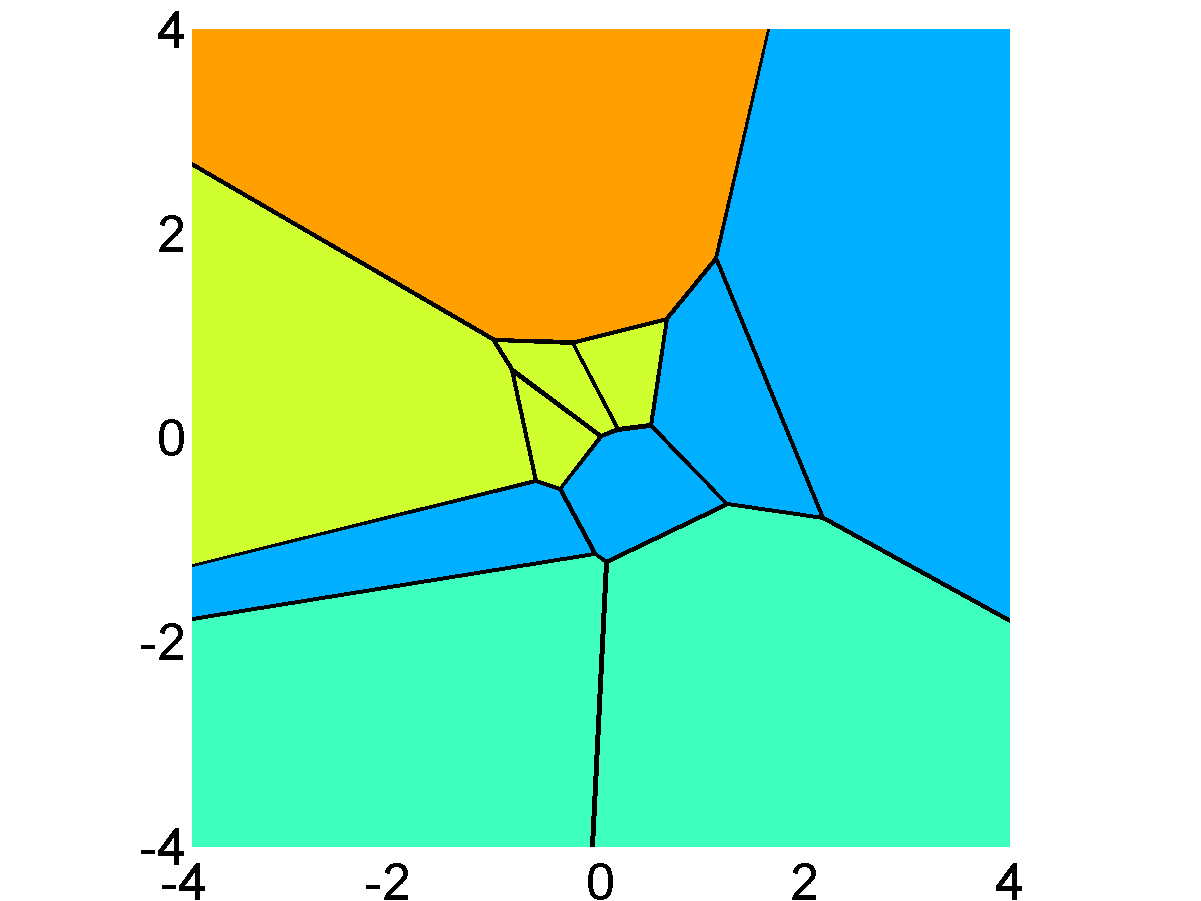} \\
&
\footnotesize  Iteration 0 &
\footnotesize  3000 &
\footnotesize  6000 &
\footnotesize  9000 \\
\end{tabular}
\caption{States of the simulated annealing
\label{fig:states_sim_anneal}
}
\end{figure}

The same covariance matrices are used as in the synthetic realization Fig.~\ref{fig:synt_all}(b,c) for the empirical transition probabilities. The states of the simulated annealing are shown at iteration 0, 3000, 6000, 9000, in Fig.~\ref{fig:states_sim_anneal} for two Markov chains with different initial states.
Despite the large difference in the initial states, the states of the truncation maps are similar after 9000 iterations.

\begin{figure}[!ht]
\begin{tabular}{ccc}
\includegraphics[width= 0.3\textwidth]{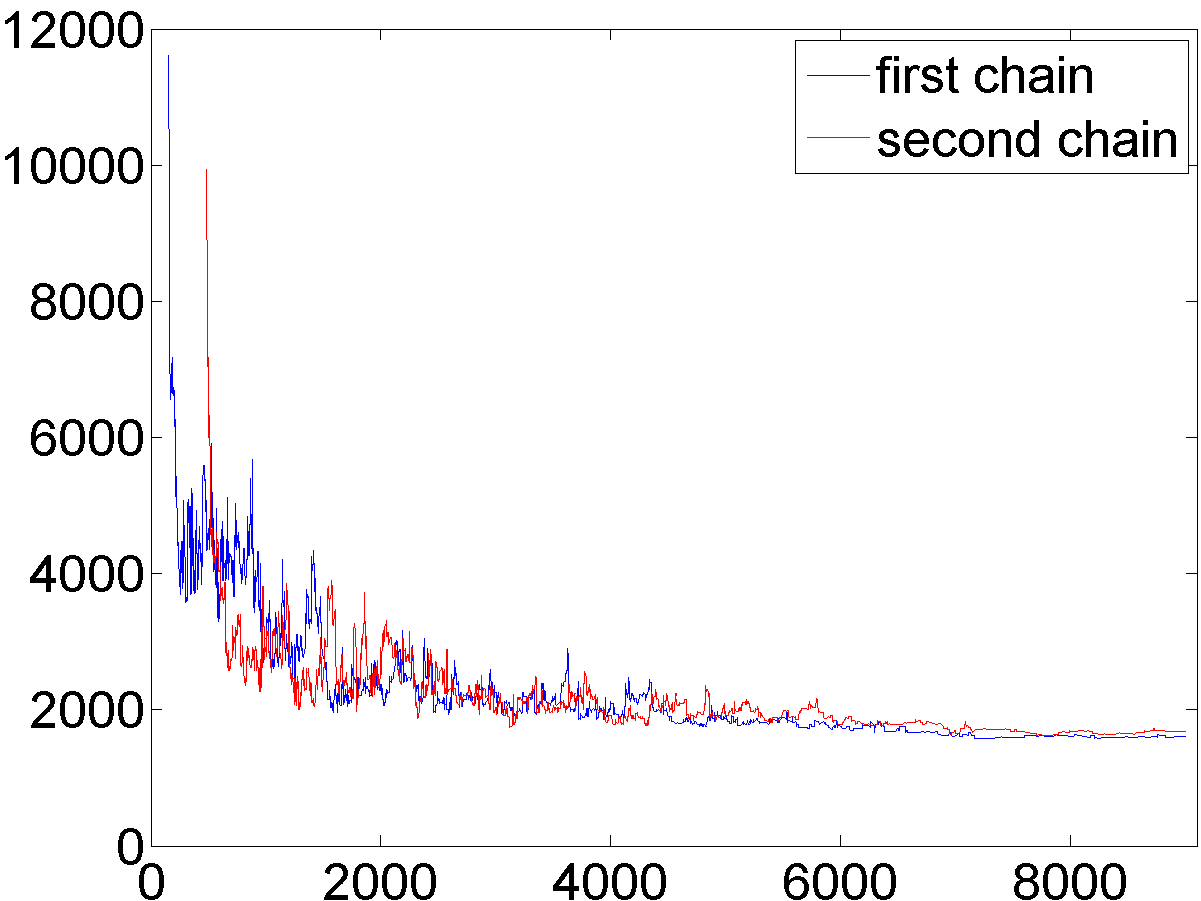} &
\includegraphics[width= 0.3\textwidth]{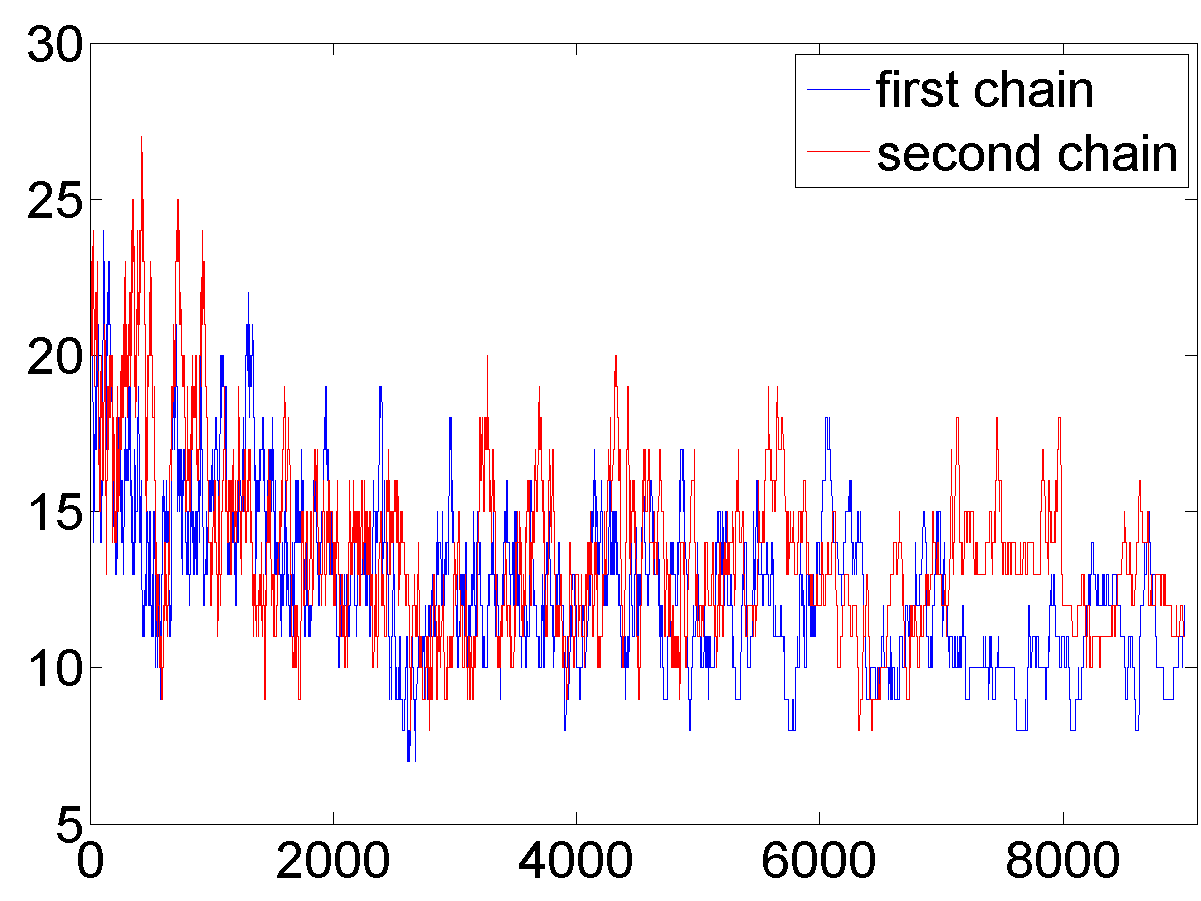} &
\includegraphics[width= 0.3\textwidth]{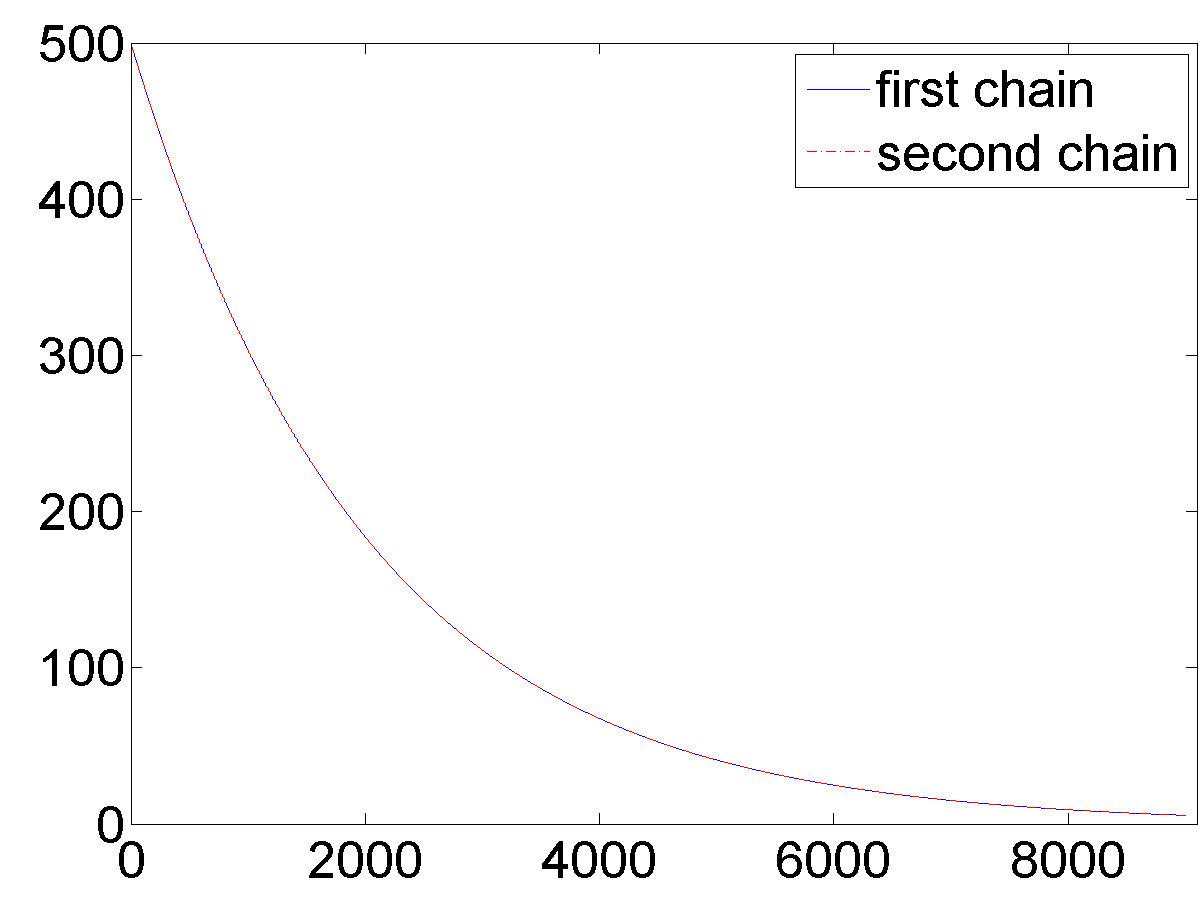}  \\
\footnotesize(a) Likelihood function $F(\theta^{(t)})$  &
\footnotesize(b) Number of nodes for $\theta^{(t)}$ &
\footnotesize(c) Temperature $T^{(t)}$  \\
\end{tabular}
\caption{Simulated annealing behavior at iteration $t=1,\dots,9000$.
\label{fig:other_sim_anneal}
}
\end{figure}

Other parameters behavior related to the estimation process are shown in Fig.~\ref{fig:other_sim_anneal}. Figures~\ref{fig:other_sim_anneal}(a,b) show the likelihood function values and the number of nodes, respectively, plotted for  iteration 1 to 9000. The temperature cooling schedule $T^{(t)} = T_0 \alpha^t, \alpha=0.9995,T_0=500$ is visualized for the same range of iterations in Fig.~\ref{fig:other_sim_anneal}(c).
The likelihood functions decline with the iterations to a cooled state with close likelihood values. The mismatch of the initial states and some first few iterations for both chains is computationally to large to be tracked. At early iterations all the proposals are accepted. The prior mean number of nodes equals 20 here, which is well reflected at the early iterations, Fig.~\ref{fig:other_sim_anneal}(b). Then the number of nodes declines, presumably due to the fact that the the generation of the number of nodes does not take into account the probabilities related to their colors. This behavior is similar for both chains.  

\subsection{Validation of results}

\begin{figure}[!ht]
\centering
\begin{tabular}{cc}
\includegraphics[width= 0.27\textwidth]{ref_grid.png} & \includegraphics[width= 0.27\textwidth]{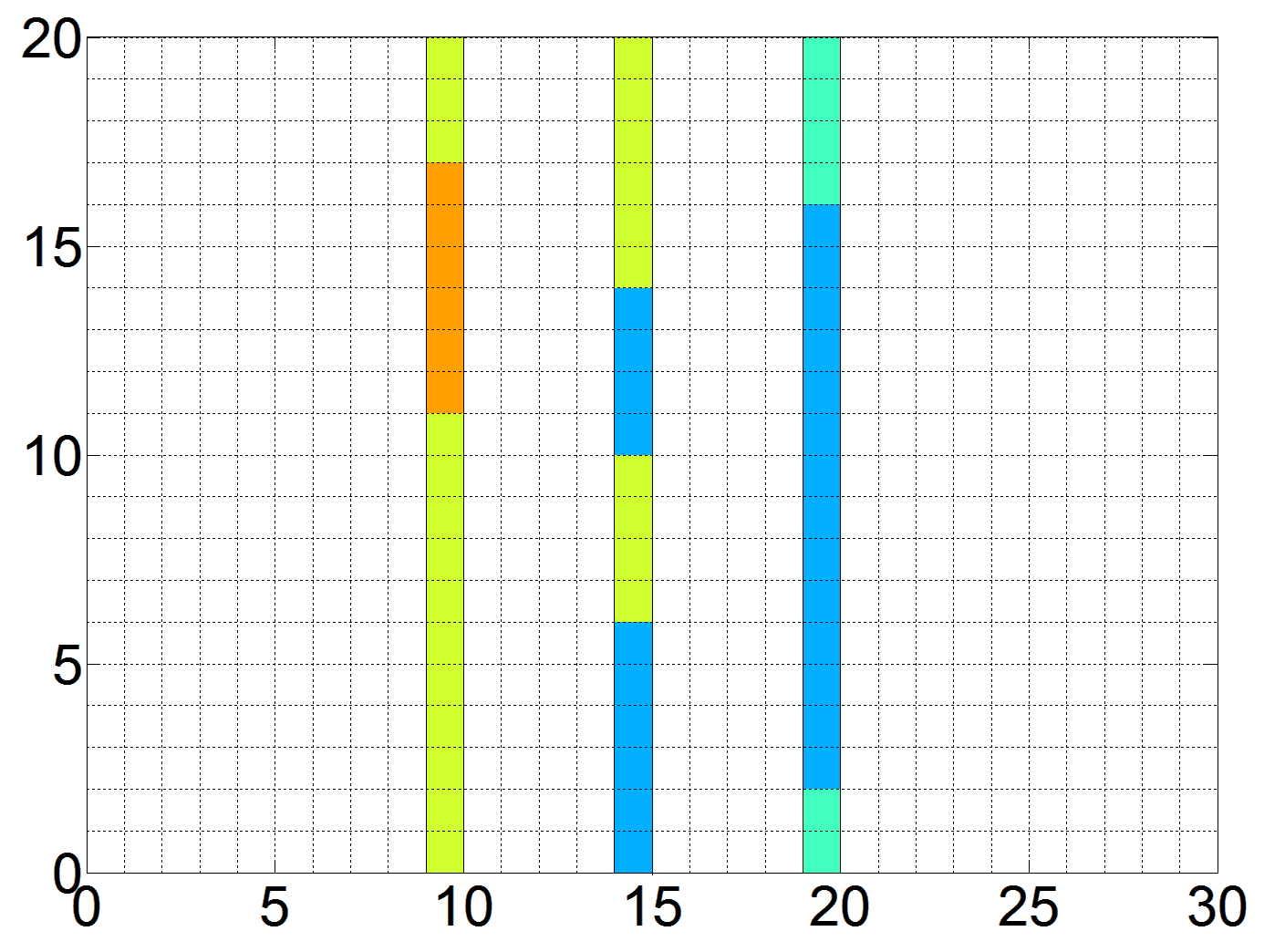} \\
{\footnotesize(a) Synthetics field} & 
{\footnotesize(b) Event} \\
\end{tabular}
\caption{Well categorical observations (b) at the cross-sections $X=10,15,20$ of the synthetics field (a)
\label{fig:ref_wells} 
}
\end{figure}

Estimation results are validated based on a set of categorical observations at the cross-sections $X=10,15,20$ of the synthetic field, Fig.~\ref{fig:ref_wells}(a), the same field  that was used to generate the transition data. It is reproduced here for visualization.
Fig.~\ref{fig:ref_wells}(b) shows the materialized event. The event includes three aligned columns of categorical observations resembling vertical well categorical data. This event is used for the validation of the predictive distributions based on the truncated bivariate model parameters estimated above. 

Because the same covariance matrices are used  in the synthetic realization for the empirical transition probabilities and in the prior model, the validation does not evaluate the latent variables distribution assumptions,
although it evaluates the results of the estimation of the truncation map based on likelihood to the transition probability data. 
The validation also answers the question of how appropriate was the estimation based on the bivariate probabilities in regard to the long correlated observation vector in this particular case.

\begin{figure}[!ht]
\begin{tabular}{ccc}
\includegraphics[trim= 0 25mm 0 0, clip=true,width= 0.3\textwidth]{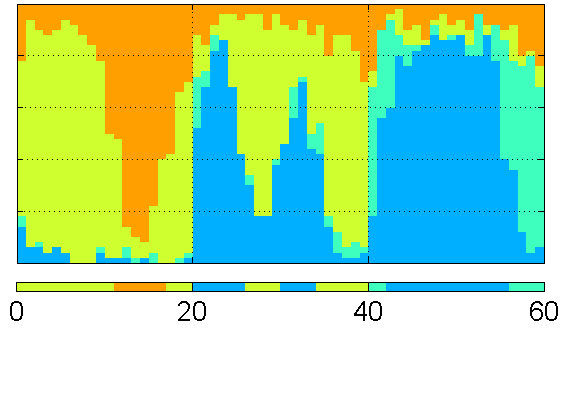} &
\includegraphics[trim= 0 25mm 0 0, clip=true,width= 0.3\textwidth]{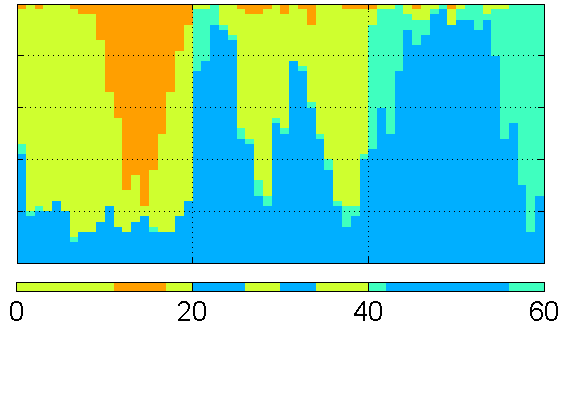}  &
\includegraphics[trim= 0 25mm 0 0, clip=true,width= 0.3\textwidth]{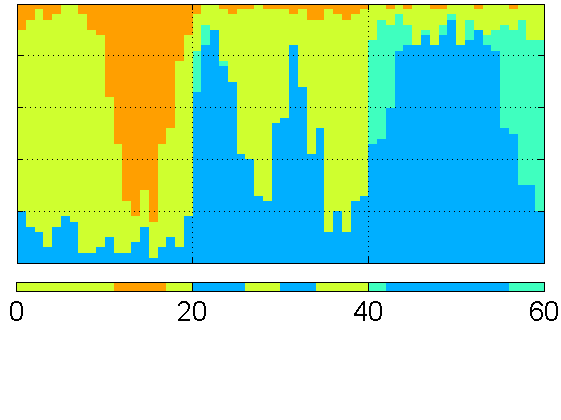}  \\
\footnotesize(a) Reference truncation map &
\footnotesize(b) First chain &
\footnotesize(b) Second chain \\
\end{tabular}
\caption{ Empirical predictive distributions based on random subsets of observations for the reference truncation map (a) and the estimated truncation maps (a,b)
\label{fig:sc_rules_ref_est}
}
\end{figure}

The predictive distributions are approximated through samples conditioned to random subsets of observations. The number of random subsets for empirical predictive distribution at each location equals 50.
An example of the simulation conditioned on subsets of the data is given in Appendix B.  
The reference model truncation map is available  in this example, Fig.~\ref{fig:synt_all}(d).
The predictive distributions based on the refernce model is shown in Fig.~\ref{fig:sc_rules_ref_est}(a).
Empirical distributions for the predictive distributions based on the estimated truncation map from each chain are given in Fig.~\ref{fig:sc_rules_ref_est}(b,c). 
The visualizations of the categorical predictive distributions from the reference map and from the estimates should be interpreted as a portion of a color in a column for each observation aligned in horizontal-axis below. 
The columns 1 to 20 correspond to the  well observations at the cross-section $X=10$, the columns 21 to 40 correspond  to the  observations at cross-section $X=15$, and the last columns from 41 to 60 are related the observations at the cross-section $X=20$. 
A relative height of the observation color (category) from the horizontal-axis occupies in the respective column above it represents the probability of a correct prediction at the respective location. Other categories have the probability proportional to their heights in this column. 

Validation of probabilistic predictions is accomplished using 
scoring rules.
If the value of the scoring rule is relatively small, the simulated 
annealing parameters can be assumed to be set correctly, for example, the 
temperature cooling schedule was slow enough. In this synthetic example, 
for which the true model parameters are known, it is possible to compare 
the scoring rule for probabilistic forecasts obtained using parameters 
from the reference model with the scoring rule computed for 
probabilistic forecasts obtained using parameters $\theta$  from the estimated 
maximum likelihood model. 
The respective scoring rules for the three predictive distributions, approximating Eq.~\eqref{eq:unordered_sc} are $S_{\text{ref}} = -23.214$;  $S_{\text{1st chain}} = -22.299$;  $S_{\text{2nd chain}} = -21.615$. The predictive power of three distribution is similar. The predictive distribution from the reference truncation map in Fig.~\ref{fig:sc_rules_ref_est}(a) overestimates the `orange' category. the predictive distribution from the estimated truncation map (first chain), Fig.~\ref{fig:sc_rules_ref_est}(b) overestimates the `blue' category. 
The predictive distribution from the estimated truncation map (second chain), Fig.~\ref{fig:sc_rules_ref_est}(c) is close to the distribution given the first chain estimate, both visually and based on the scoring rule values.
Overall, the validation shows that the estimated truncation maps are comparable to the reference truncation map, despite the small number of samples. 

\section{Conclusions}

One of the challenges of using the TPG is the problem of estimating parameters of the truncation map that will generate realizations of categorical vectors with the desired distribution. 
This paper proposed an approach for estimation of a truncation map  using a Voronoi tessellation to define regions of the truncation map corresponding to different categorical variables. The categorical bivariate unit-lag distribution   was assumed to be known from geological data. Estimation of parameters of the truncation map was performed using simulated annealing to minimize an objective function measuring the difference between a joint distributions from the BME model used to integrate the unit-lag transition probabilities in different directions, and the respective empirical distribution obtained from sampling with the truncation map. The  covariance matrices of the  latent GRFs were assumed to be known during the optimization. Parameters that were estimated   included the number of Voronoi nodes, their locations and the categories assigned to each node. 
The validation of the result was performed by the scoring rules computation based on the unordered data from synthetic categorical joint observation as an event.
An efficient method of conditional simulation for the defined model, was based on the propagative version of the Gibbs sampler. In particular, the simulation method was used in the validation implementation for conditioning to a number of observations with strong correlations, and exhibited good mixing properties.
The validation of the estimates of the truncation map based on the scoring rules has been successful with the predictive power of the estimated model similar to each other and to the truncation map used for creating the synthetic field.

\section{Acknowledgments}
The authors of this paper would like to thank the IRIS / Uni Research CIPR cooperative research project 
`Reservoir Data Assimilation For Realistic Geology' funded by its industry partners ConocoPhillips, Eni, Petrobras, Statoil, and Total, as well as the Research Council of Norway (PETROMAKS).

\noindent
\appendix

\numberwithin{equation}{section}
\numberwithin{figure}{section}
\numberwithin{algorithm}{section}

\renewcommand{\thesection}{Appendix \Alph{section}:}
\section{Proof of Csisz\'ar criterion} 
\renewcommand{\thesection}{\Alph{section}}

A proof of the first equality:
\begin{align*} 
D ( q  \|  p^* ) &= \sum_z q(z) \ln \frac{q(z)}{p^*(z)} 
= \sum_{z} q(z_B) q(z_{B_*\backslash B}|z_B) \ln \frac{q(z)/\pi(z_B)}{p^*(z)/\pi(z_B)} \\
&= \sum_{z_B} q(z_B) \sum_{B_*\backslash B} q(z_{B_*\backslash B}|z_B) \ln \frac{q(z_{B_*\backslash B}|z_B)}{z^*(z_{B_*\backslash B}|z_B)} \\
&=\sum_{z_B} q (z_B) D \bigl( q ( \cdot_{B_*\backslash B} | \cdot_B ) \|  p ( \cdot_{B_*\backslash B} | \cdot_B ) \bigr).
\end{align*}
A proof of the second equality:
\begin{align*} 
D ( p^* \|  p ) &= \sum_{z} p^* (z) \ln \dfrac{p^* (z)}{p (z)} 
= \sum_{z} \pi (z_B) p (z_{B_*\backslash B} | z_B ) \ln 
\dfrac{ \pi (z_B) }{ p (z_B) } \\
&= \sum_{z_B} \pi (z_B) \ln \dfrac{ \pi (z_B) }{ p (z_B) } = 
D \bigl( p^* (\cdot_B ) \|  p (\cdot_B) \bigr). 
\end{align*}

\renewcommand{\thesection}{Appendix \Alph{section}:}
\section{Conditional simulation example} 
\renewcommand{\thesection}{\Alph{section}}

\label{sec:cond_sim_ex}
\begin{figure}[!ht]
\begin{tabular}{ccc}
\includegraphics[width= 0.3\textwidth]{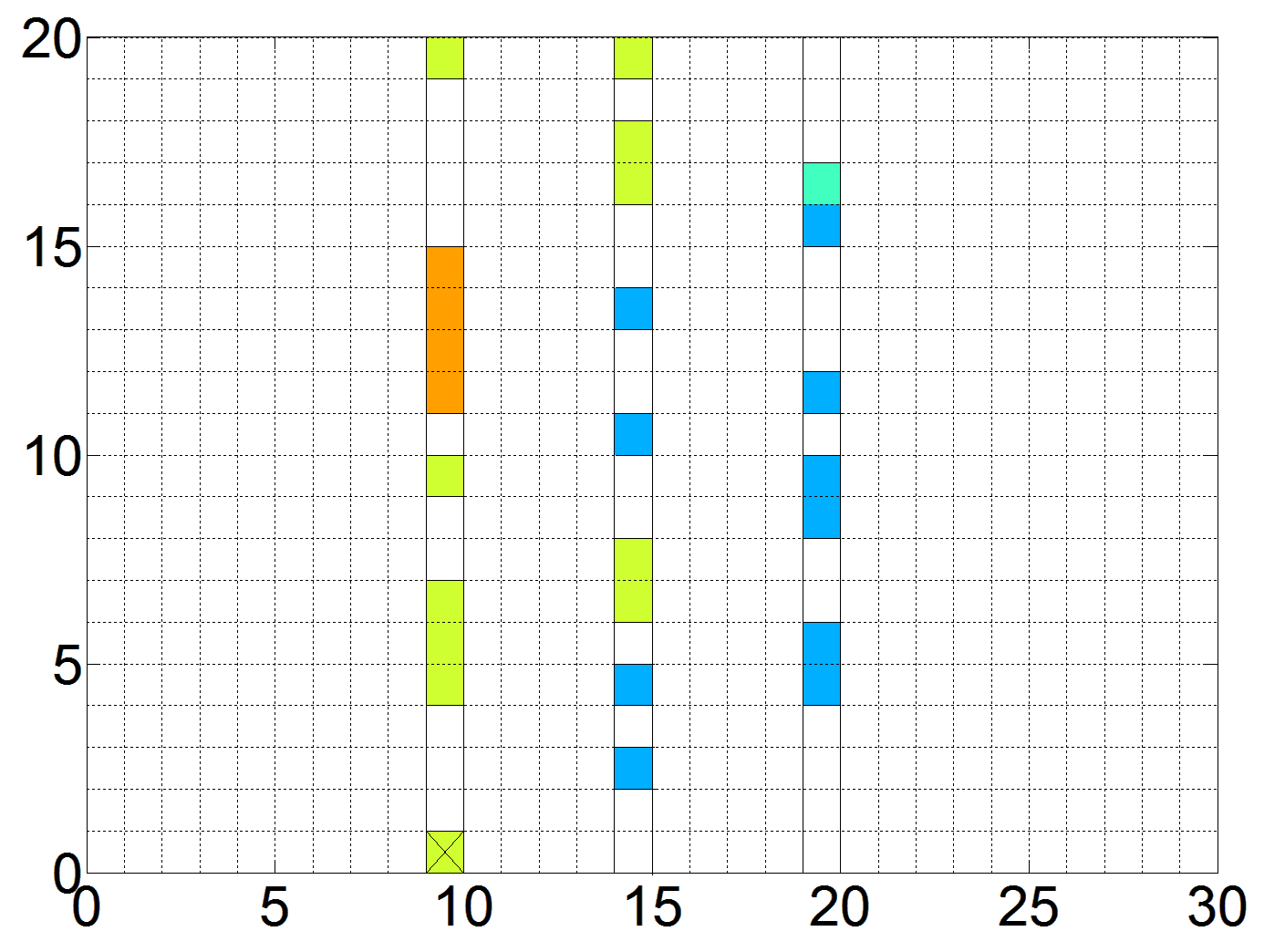} &
\includegraphics[width= 0.3\textwidth]{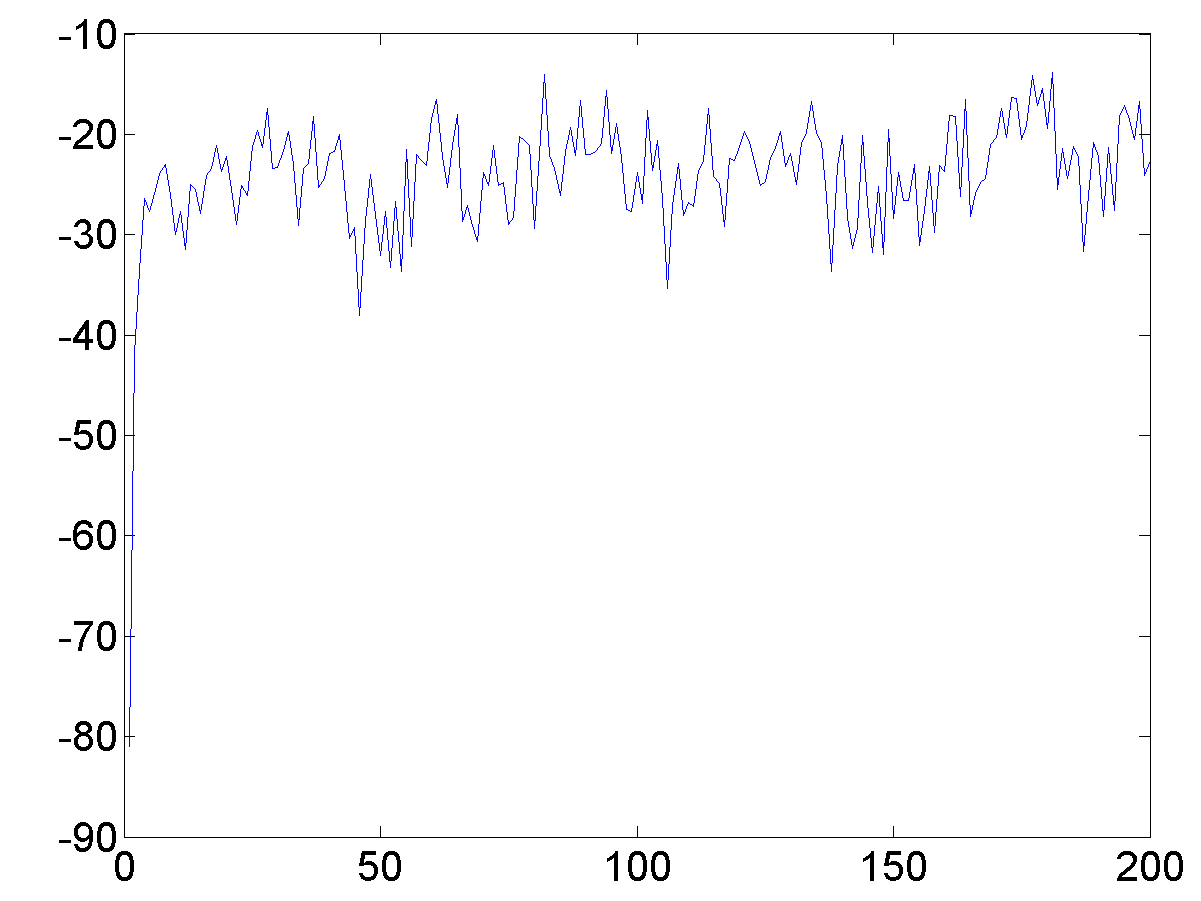} &
\includegraphics[width= 0.3\textwidth]{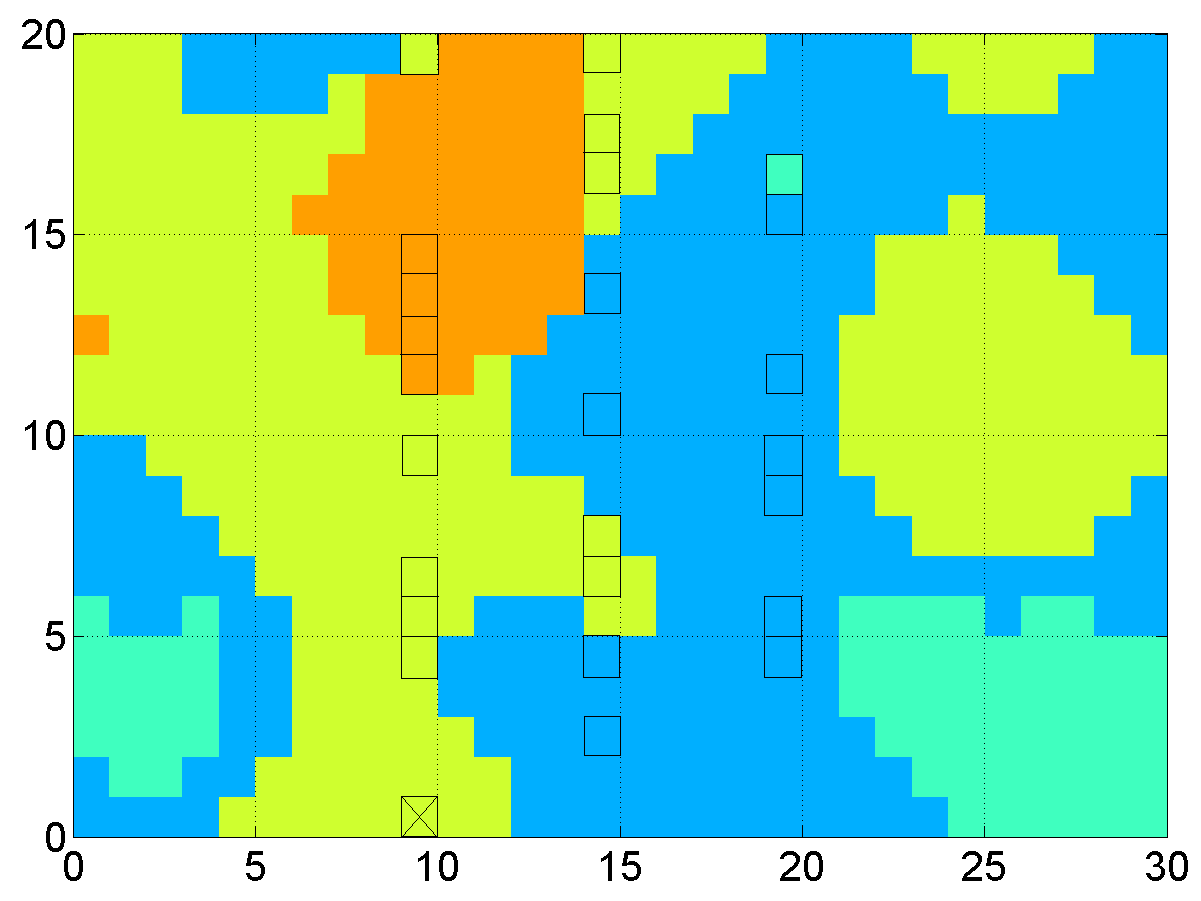} \\
\footnotesize(a) &
\footnotesize(b) & 
\footnotesize(c)  \\
\end{tabular}
\caption{ (a) Conditioning to the random subset of observations, (b) mixing of the latent Gaussian vector realizations during 200 iterations, (c) conditional simulation on the grid
\label{fig:unordered_cond} 
}
\end{figure}

This part provides an example of conditional simulation based on a random subset of the observations in the event, Fig.~\ref{fig:ref_wells}, and the estimated truncation map Fig.~\ref{fig:states_sim_anneal} (first chain). To approximate a member of the sum for the scoring rule in Eq.~\eqref{eq:unordered_sc}, multiple random subsets should be sampled, and the Gaussian vectors, Eq.~\eqref{eq:cond_sim} should be simulated as one of the steps. An example is related to the member of the sum for the first observation of the event with the coordinates $(10,1)$. Figure~\ref{fig:unordered_cond}(a) shows the unordered subset of the observations. The latent Gaussian vectors with the related coordinates and the covariance matrix, is conditioned to these categorical observations according to the method described. 200 iterations of both the standard and alternative Gibbs sampler were used to ensure that the Markov chain for the Gaussian realizations have mixed Fig.~\ref{fig:unordered_cond}(b). The vertical axes represents the sum of the logarithms of the  probability densities of the respective latent Gaussian random vector. 
After the conditional Gaussian realizations at the related coordinates are obtained, the Gaussian simulation allows to sample and map the categorical realizations at the coordinate $(10,1)$ (crossed square). For the sake of the example, the sampling is done in the entire grid, Fig.~\ref{fig:unordered_cond}(c). In this realization example the simulated categorical observation at $(10,1)$ is predicted correctly.


\begin{thebibliography}{}

\bibitem[\protect\astroncite{Al-Anezi et~al.}{2013}]{alanezi:13}
Al-Anezi, K., Kumar, S., Ebaid, A., Bonnel, A., Lucet, N., Lecante, G., and
  Ortet, S. (2013).
\newblock Geostatistical modeling with seismic characterization of
  {W}ara/{B}urgan sands {M}inagish {F}ield {W}est {K}uwait.
\newblock In {\em SPE Reserv. Charact. and Simul. Conf. and Exhib., Abu Dhabi,
  16--18 Sept.}, page SPE 166046. Society of Petroleum Engineers.

\bibitem[\protect\astroncite{Albert{\~a}o et~al.}{2005}]{albertao:05}
Albert{\~a}o, G.~A., Grell, A.~P., Badolato, D., and dos Santos, L.~R. (2005).
\newblock 3{D} geological modeling in a turbidite system with complex
  stratigraphic-structural framework---an example from {C}ampos {B}asin
  {B}razil.
\newblock In {\em SPE Annu. Tech. Conf. and Exhib., Dallas, Texas, 9--12 Oct.},
  page SPE 95612. Society of Petroleum Engineers.

\bibitem[\protect\astroncite{Allard et~al.}{2012}]{allard:12}
Allard, D., D'Or, D., Biver, P., and Froidevaux, R. (2012).
\newblock Non-parametric diagrams for pluri-gaussian simulations of
  lithologies.
\newblock In {\em 9th {I}nt. {G}eostat. {C}ongr., Oslo, Norway}, pages 11--15.

\bibitem[\protect\astroncite{Armstrong et~al.}{2011}]{armstrong:11}
Armstrong, M., Galli, A., Beucher, H., {Le Loc'h}, G., Renard, D., Doligez, B.,
  Eschard, R., and Geffroy, F. (2011).
\newblock {\em Plurigaussian Simulations in Geosciences}.
\newblock Springer Berlin Heidelberg, 2nd revised edition edition.

\bibitem[\protect\astroncite{Bregman}{1967}]{bregman:67}
Bregman, L.~M. (1967).
\newblock The relaxation method of finding the common point of convex sets and
  its application to the solution of problems in convex programming.
\newblock {\em USSR Comput. Math. and Math. Phys.}, 7(3):200--217.

\bibitem[\protect\astroncite{Carrillat et~al.}{2010}]{carrillat:10}
Carrillat, A., Sharma, S.~K., Grossmann, T., Iskenova, G., Friedel, T., et~al.
  (2010).
\newblock Geomodeling of giant carbonate oilfields with a new multipoint
  statistics workflow.
\newblock In {\em Abu Dhabi Int. Petroleum Exhib. and Conf.}, page
  {SPE}~137958. Society of Petroleum Engineers.

\bibitem[\protect\astroncite{Cherbunin et~al.}{2009}]{cherubini:09}
Cherbunin, C., Giasi, C.~I., Musci, F., and Pastore, N. (2009).
\newblock Application of truncated plurigaussian method for the reactive
  transport modeling of a contaminated aquifer.
\newblock In {\em Proc. of the 4th {IASME} / {WSEAS} {I}nt. {C}onf. on {W}ater
  {R}esour., {H}ydraul. \& {H}ydrol. ({WHH}'09)}, pages 119--124.

\bibitem[\protect\astroncite{Csisz{\'a}r}{1975}]{csiszar:75}
Csisz{\'a}r, I. (1975).
\newblock I-divergence geometry of probability distributions and minimization
  problems.
\newblock {\em The Ann. of Probab.}, pages 146--158.

\bibitem[\protect\astroncite{Deming and Stephan}{1940}]{deming:40}
Deming, W.~E. and Stephan, F.~F. (1940).
\newblock On a least squares adjustment of a sampled frequency table when the
  expected marginal totals are known.
\newblock {\em The Annals of Mathematical Statistics}, 11(4):427--444.

\bibitem[\protect\astroncite{Deutsch and Deutsch}{2014}]{deutsch:14}
Deutsch, J.~L. and Deutsch, C.~V. (2014).
\newblock A multidimensional scaling approach to enforce reproduction of
  transition probabilities in truncated plurigaussian simulation.
\newblock {\em Stoch. Environ. Res. and Risk Assess.}, 28(3):707--716.

\bibitem[\protect\astroncite{Du et~al.}{1999}]{du:99}
Du, Q., Faber, V., and Gunzburger, M. (1999).
\newblock Centroidal {V}oronoi tessellations: Applications and algorithms.
\newblock {\em SIAM Rev.}, 41(4):637--676.

\bibitem[\protect\astroncite{Emery}{2010}]{emery:10}
Emery, X. (2010).
\newblock On the existence of mosaic and indicator random fields with
  spherical, circular, and triangular variograms.
\newblock {\em Math. Geosci.}, 42:969--984.

\bibitem[\protect\astroncite{Emery et~al.}{2014}]{emery:14}
Emery, X., Arroyo, D., and Pel{\'a}ez, M. (2014).
\newblock Simulating large {G}aussian random vectors subject to inequality
  constraints by {G}ibbs sampling.
\newblock {\em Math. Geosci.}, 46(3):265--283.

\bibitem[\protect\astroncite{Fachri et~al.}{2013}]{fachri:13}
Fachri, M., Tveranger, J., Braathen, A., and Schueller, S. (2013).
\newblock Sensitivity of fluid flow to deformation-band damage zone
  heterogeneities: A study using fault facies and truncated gaussian
  simulation.
\newblock {\em J. of Struct. Geol.}, 52:60--79.

\bibitem[\protect\astroncite{Galli et~al.}{1994}]{galli:94}
Galli, A., Beucher, H., {Le Loc'h}, G., Doligez, B., and Group, H. (1994).
\newblock The pros and cons of the truncated {Gaussian} method.
\newblock In {\em Geostat. Simul.}, pages 217--233. Kluwer Academic, Dordrecht.

\bibitem[\protect\astroncite{Galli et~al.}{2006}]{galli:06}
Galli, A., {Le Loc'h}, G., Geffroy, F., and Eschard, R. (2006).
\newblock An application of the truncated pluri-gaussian method for modeling
  geology.
\newblock In Coburn, T.~C., Yarus, J.~M., and Chambers, R.~I., editors, {\em
  Stoch. Model. and Geostat.: Princ., Methods, and Case Stud., Vol. II: AAPG
  Comput. Appl. in Geol.}, pages 109--122. AAPG Special Volumes.

\bibitem[\protect\astroncite{Gneiting and Raftery}{2007}]{gneiting:07a}
Gneiting, T. and Raftery, A.~E. (2007).
\newblock Strictly proper scoring rules, prediction, and estimation.
\newblock {\em Journal of the American Statistical Association},
  102(477):359--378.

\bibitem[\protect\astroncite{Jaynes}{1957}]{jaynes:57}
Jaynes, E.~T. (1957).
\newblock Information theory and statistical mechanics.
\newblock {\em Phys. Rev.}, 106(4):620--630.

\bibitem[\protect\astroncite{Ki{\^e}u et~al.}{2013}]{kieu:13}
Ki{\^e}u, K., Adamczyk-Chauvat, K., Monod, H., and Stoica, R.~S. (2013).
\newblock A completely random {T}-tessellation model and {G}ibbsian extensions.
\newblock {\em Spatial Stat.}, 6(0):118 -- 138.

\bibitem[\protect\astroncite{Kyriakidis et~al.}{1999}]{kyriakidis:99}
Kyriakidis, P.~C., Deutsch, C.~V., and Grant, M.~L. (1999).
\newblock Calculation of the normal scores variogram used for truncated
  {G}aussian lithofacies simulation: theory and {FORTRAN} code.
\newblock {\em Comput. \& Geosci.}, 25(2):161--169.

\bibitem[\protect\astroncite{Lantu{\'e}joul and Desassis}{2012}]{lantuejoul:12}
Lantu{\'e}joul, C. and Desassis, N. (2012).
\newblock {Simulation of a Gaussian random vector: a propagative approach to
  the Gibbs sampler}.
\newblock In {\em {9th Int. Geostat. Congr.}}, Oslo, Norway.

\bibitem[\protect\astroncite{{Le Loc'h} et~al.}{1994}]{leloch:94}
{Le Loc'h}, G., Beucher, H., Galli, A., Doligez, B., and Group, H. (1994).
\newblock Improvement in the truncated {Gaussian} method: Combining several
  {Gaussian Functions}.
\newblock In {\em Proc. of the 4th Eur. Conf. on the Math. of Oil Recovery
  (ECMOR)}, page 13 p.

\bibitem[\protect\astroncite{Mariethoz et~al.}{2009}]{mariethoz:09}
Mariethoz, G., Renard, P., Cornaton, F., and Jaquet, O. (2009).
\newblock Truncated plurigaussian simulations to characterize aquifer
  heterogeneity.
\newblock {\em Ground Water}, 47(1):13--24.

\bibitem[\protect\astroncite{Matheron et~al.}{1987}]{matheron:87}
Matheron, G., Beucher, H., de~Fouquet, C., Galli, A., Guerillot, D., and
  Ravenne, C. (1987).
\newblock Conditional simulation of the geometry of fluvio-deltaic reservoirs,
  {SPE} 16753.
\newblock In {\em Proc. of the 62nd Annu. Tech. Conf. of the SPE}, pages
  123--131.

\bibitem[\protect\astroncite{Perulero~Serrano et~al.}{2012}]{perulero:12}
Perulero~Serrano, R., Guadagnini, L., Giudici, M., Guadagnini, A., and Riva, M.
  (2012).
\newblock Application of the truncated plurigaussian method to delineate
  hydrofacies distribution in heterogeneous aquifers.
\newblock In {\em The 19th {I}nt. {C}onf. on {W}ater {R}esour. ({CMWR})}, pages
  952--959.

\bibitem[\protect\astroncite{Perulero~Serrano et~al.}{2014}]{perulero:14}
Perulero~Serrano, R., Guadagnini, L., Riva, M., Giudici, M., and Guadagnini, A.
  (2014).
\newblock Impact of two geostatistical hydro-facies simulation strategies on
  head statistics under non-uniform groundwater flow.
\newblock {\em J. of Hydrol.}, 508:343--355.

\bibitem[\protect\astroncite{Tarantola}{2005}]{tarantola:05}
Tarantola, A. (2005).
\newblock {\em Inverse Problem Theory and Methods for Model Parameter
  Estimation}.
\newblock Society for Industrial Mathematics.

\bibitem[\protect\astroncite{Xu et~al.}{2006}]{xu:06}
Xu, C., Dowd, P.~A., Mardia, K.~V., and Fowell, R.~J. (2006).
\newblock A flexible true plurigaussian code for spatial facies simulations.
\newblock {\em Comput. \& Geosci.}, 32(10):1629--1645.

\bibitem[\protect\astroncite{Zelen and Severo}{2012}]{zelen:12}
Zelen, M. and Severo, N.~C. (2012).
\newblock Probability functions.
\newblock In Abramowitz, M. and Stegun, I.~A., editors, {\em Handb. of Math.
  Funct.: with Formulas, Graphs, and Math. Tables}, pages 925--995. Courier
  Dover Publications.

\end{thebibliography}
\end{document}